\documentclass[onefignum,onetabnum]{siamart220329}
\usepackage{lipsum}
\usepackage{amssymb,amsfonts}
\usepackage{amsmath}
\usepackage{graphicx,subfigure}
\usepackage{epstopdf}
\usepackage{algorithm,algorithmic}
\usepackage{color}
\usepackage{bm}
\usepackage{booktabs,array}
\usepackage{multirow,multicol}

\usepackage[contents = Peer~review~only, scale = 8, color = black, angle = 50, opacity = .10]{background}
\ifpdf
\DeclareGraphicsExtensions{.eps,.pdf,.png,.jpg}
\else
\DeclareGraphicsExtensions{.eps}
\fi

\newsiamremark{remark}{Remark}
\newsiamremark{hypothesis}{Hypothesis}
\crefname{hypothesis}{Hypothesis}{Hypotheses}
\newsiamthm{claim}{Claim}

\headers{DLRA for strongly damped wave equations}{Y.-L. Zhao and X.-M. Gu}

\title{A low-rank algorithm for strongly damped wave equations with visco-elastic damping and mass terms\thanks{Received by the editors;
		\funding{This work was supported by the Natural Science Foundation of Sichuan Province (No.~2022NSFSC1815)
			and the Sichuan Science and Technology Program (No.~2023NSFSC1326).}}}

\author{Yong-Liang Zhao\thanks{School of Mathematical Sciences, 
		Sichuan Normal University, Chengdu, Sichuan 610066, P.R. China 
		(\email{ylzhaofde@sina.com}).}
	\and Xian-Ming Gu\thanks{School of Economic Mathematics,
		Southwestern University of Finance and Economics,
		Chengdu, Sichuan 611130, P.R. China
	(\email{guxianming@live.cn}).}
}

\usepackage{amsopn}

\ifpdf
\hypersetup{
  pdftitle={DLRA for strongly damped wave equations},
  pdfauthor={Y.-L. Zhao}
}
\fi

\begin{document}

\maketitle

\begin{abstract}
Damped wave equations have been used in many real-world fields. 
In this paper, we study a low-rank solution of the strongly damped wave equation with the damping term, 
visco-elastic damping term and mass term.
Firstly, a second-order finite difference method is employed for spatial discretization.
Then, we receive a second-order matrix differential system.
Next, we transform it into an equivalent first-order matrix differential system,
and split the transformed system into three subproblems.
Applying a Strang splitting to these subproblems and combining a dynamical low-rank approach,
we obtain a low-rank algorithm.
Numerical experiments are reported to demonstrate that the proposed low-rank algorithm 
is robust and accurate, and has second-order convergence rate in time.
\end{abstract}



\section{Introduction}
\label{sec1}

Wave equations with damping terms have been used in many fields 
such as visco-elastic materials \cite{pata2006smooth}, 
quantum mechanics \cite{ghidaglia1991longtime}, railway tracks \cite{ansari2011frequency,edalatzadeh2018stability}
and fluid dynamics \cite{ponce1985global}.
Now, damped wave equations received many researchers' attention for their rich physical significance.
Numerous studies about the properties of these equations have been carried out, 
e.g., the existence and blow-up of solutions, attractors and long-time behavior,
see \cite{ghidaglia1991longtime,ponce1985global,GAZZOLA2006185,LianXu+2020+613+632,pata2005strongly,
	carvalho2008strongly,xie2009asymptotic,kalantarov2009finite,yang2010exponential,lian2019global,
	chen2021blow,wang2023long,ammari2023uniform} and references therein. 
For instance, Ponce \cite{ponce1985global} studied the global existence of solutions 
of some quasilinear equations 
such as the dissipative equation with linear part $u_{tt} - \Delta u - v u_t$.
For strongly damped abstract nonlinear wave equations, 
Ghidaglia and Marzocchi \cite{ghidaglia1991longtime} gave uniform time estimates 
and constructed global attractors.
For the strongly damped semilinear wave equation, 
Pata and Squassina \cite{pata2005strongly} proved the existence of the universal attractor.
They also revealed that the asymptotic behavior of its solutions 
in connection with the damping coefficient.
Gazzola and Squassina \cite{GAZZOLA2006185} showed that for damped wave equations with superlinear term,
every global solution is uniformly bounded in the phase-space.
For arbitrarily high energy initial data, they proved that some solutions of the considered model
blow up in finite time.
For strongly damped wave equations with power nonlinear term 
(i.e., $|u_t|^q$ or $|u|^p + |u_t|^q$, $p,q > 1$) in an exterior domain,
under some assumptions on the initial value and the exponents $p,q$, 
Chen and Fino \cite{chen2021blow} proved: i) the local existence of the mild solutions; 
ii) the blowing up of the solution in finite time.

Generally, it is hard to reach the analytical solutions of (strongly) damped wave equations 
due to their mathematical complexity.
Thus, numerical solutions of such equations attract much attention.
Larsson et al.~\cite{larsson1991finite} proposed a class of finite element schemes for
a strongly damped linear wave equation.
The optimal convergence orders of these schemes are also investigated.
Based on a two-step, one-parameter method, Djidjeli et al.~\cite{djidjeli1995numerical}
developed a numerical method for a two-dimensional (2D) perturbed Sine-Gordon equation.
For the equation considered in \cite{larsson1991finite}, 
Thom{\'e}e and Wahlbin \cite{thomee2004maximum} derived two fully discrete schemes,
which are established by employing piecewise linear finite elements in space
and backward Euler and Crank–Nicolson methods in time. 
The optimal order error estimations in the maximum-norm of their schemes are conducted.
Achouri \cite{achouri2019finite} proposed two second-order finite difference schemes for
the 2D damped semilinear wave equation.
The uniqueness of solutions, unconditional stability and convergence 
of these two schemes are rigorously analyzed.
Combining a generalized finite element (FE) method and a backward Euler method,
Ljung et al.~\cite{ljung2021generalized} designed a FE scheme for 
the strongly damped wave equation with highly varying coefficients.
Other numerical methods for (strongly) damped wave equations have been studied, see
\cite{wang2005numerical,achouri2021analysis,phan2022exponential}.

Dynamical low-rank (DLR) approximations as a model reduction technique 
have a wide range of application areas, 
see \cite{koch2007dynamical,lubich2014projector,Ceruti2022unconventional,Ceruti2022adaptiveDLR,
	einkemmer2023robust,carrel2023low,ostermann2019convergence,zhao2021low} and references therein
for more about it. 
Recently, Hochbruck et al.~\cite{hochbruck2023dynamical,hochbruck2023rank} conducted 
new DLR integrators for second-order matrix differential equations.
Inspired by these works,
in this paper, we design a low-rank approximation for the following 
strongly damped wave equation \cite{phan2022exponential} 
describing the behavior of waves that lose energy over time:
\begin{equation}\label{eq1.1}
	\begin{cases}
		\begin{split}
			u_{tt}(x,y,t) + & \gamma u_t(x,y,t) + \delta u(x,y,t) = 
			\Delta \left[ \alpha u(x,y,t) + \beta u_t(x,y,t) \right] + \\
		    & f\left( u(x,y,t) \right) + g \left( u_t(x,y,t) \right), \qquad (x,y,t) \in \Omega \times (0,T),
		\end{split} \\
		u(x,y,0) = p(x,y), \quad u_t(x,y,0) = q(x,y), \hspace{5mm} (x,y) \in \Omega, \\
		u(x,y,t) = 0, \hspace{48mm} (x,y,t) \in \partial\Omega \times (0,T),
	\end{cases}
\end{equation}
where $\alpha > 0$, $\beta,\gamma,\delta \geq 0$, 
$\Omega = [x_L,x_R] \times [y_L,y_R] \subset \mathbb{R}^2$
and $p(x,y), q(x,y), f, g$ are known functions.
In Eq.~\cref{eq1.1}, $\gamma u_t(x,y,t)$, $\beta \Delta u_t(x,y,t)$ and $\delta u(x,y,t)$
are the damping term, the structural (visco-elastic) damping term and the mass term, respectively.
As far as we are aware, this is the first attempt in the literature to find a low-rank solution to Eq.~\cref{eq1.1}.

The rest of this paper is organized as follows. 
In \cref{sec2}, we derive the semi-discrete scheme for Eq.~\cref{eq1.1}.
It is a second-order ordinary matrix differential equation.
Then, this scheme is transformed into a first-order ordinary matrix differential equation.
\Cref{sec3} derives our low-rank algorithm by 
applying a Strang splitting and a DLR approximation to the
first-order ordinary differential equation.
In \cref{sec4}, three numerical examples are reported to
show the efficiency and accuracy of the proposed low-rank algorithm.
Some conclusions are drawn in \cref{sec5}.

\section{The ordinary matrix differential equation}
\label{sec2}

In this section, we use a finite difference method for the discretization in the spatial variable
of Eq.~\cref{eq1.1}. Then, the ordinary matrix differential equation is derived.

Let $h_x = \frac{x_R - x_L}{N_x}$ and $h_y = \frac{y_R - y_L}{N_y}$ be the grid spacing in $x$ and
$y$ directions, respectively. Here, $N_x$ and $N_y$ are two given positive integer numbers.
Then, the considered space domain $\Omega$ is discretized by 
\begin{equation*}
	\Omega_s = \left\{ (x_i,y_j) \mid x_i = x_L + i h_x, y_j = y_L + j h_y,
	0 \leq i \leq N_x, 0 \leq j \leq N_y \right\}.
\end{equation*}
For approximating the Laplacian operator, we choose a second-order finite difference method.
That is,
\begin{align*}
	\Delta u(x_i,y_j,t) & \approx \frac{u_{i - 1,j}(t) - 2 u_{ij}(t) + u_{i + 1,j}(t)}{h_x^2} + 
	\frac{u_{i,j - 1}(t) - 2 u_{ij}(t) + u_{i,j + 1}(t)}{h_y^2} \\
	& \triangleq \delta_h^2 u_{ij}(t),
\end{align*}
where $u_{ij}(t)$ is the numerical approximation of $u(x_i,y_j,t)$.
Substituting this approximation into Eq.~\eqref{eq1.1}, we get the following semi-discrete scheme:
\begin{equation*}
	\frac{d^2 u_{ij}(t)}{d t^2} + \gamma \frac{d u_{ij}(t)}{d t} + \delta u_{ij}(t) = 
	\delta_h^2 \left( \alpha u_{ij}(t) + \beta \frac{d u_{ij}(t)}{d t} \right) +
	f\left( u_{ij}(t) \right) + g\left( \frac{d u_{ij}(t)}{d t} \right).
\end{equation*}
The numerical solution of Eq.~\cref{eq1.1} can be reached 
by applying some exponential integrators \cite{phan2022exponential} 
or the leapfrog scheme (also known as the St\"{o}mer scheme or the Verlet scheme) 
\cite{hairer2003geometric,carle2022error} to the above semi-discrete scheme.
Denote 
\begin{equation*}
	P(t) = \left[ u_{i j}(t) \right]_{\substack{1 \leq i \leq N_x - 1 \\ 1 \leq j \leq N_y - 1}}, ~
	{P}'(t) = \left[ \frac{d u_{i j}(t)}{d t} \right]_{\substack{1 \leq i \leq N_x - 1 \\ 
			1 \leq j \leq N_y - 1}}, ~
	{P}''(t) = \left[ \frac{d^2 u_{i j}(t)}{d t^2} \right]_{\substack{1 \leq i \leq N_x - 1 \\ 
			1 \leq j \leq N_y - 1}}, 
\end{equation*}
\begin{equation*}
	P^0 = \left[ p(x_i,y_j) \right]_{\substack{1 \leq i \leq N_x - 1 \\ 1 \leq j \leq N_y - 1}}, \qquad
	Q^0 = \left[ q(x_i,y_j) \right]_{\substack{1 \leq i \leq N_x - 1 \\ 1 \leq j \leq N_y - 1}}.
\end{equation*}
Then, the semi-discrete scheme can be rewritten as
\begin{equation}\label{eq2.1}
	\begin{split}
	{P}''(t) + \gamma {P}'(t) + \delta P(t) = & -\alpha \left( A_x P(t) + P(t) A_y \right) - \\ 
	& \beta \left( A_x P'(t) + P'(t) A_y \right) + F\left( P(t),  {P}'(t)\right)
	\end{split}	
\end{equation} 
with initial values $P(0) = P^0$ and ${P}'(0) = Q^0$.
Here, 
\begin{align*}
	A_x = \frac{1}{h_x^2} \mathrm{tridiag}(-1,2,-1) \in \mathbb{R}^{(N_x - 1) \times (N_x - 1)}, \\
	A_y = \frac{1}{h_y^2} \mathrm{tridiag}(-1,2,-1) \in \mathbb{R}^{(N_y - 1) \times (N_y - 1)}
\end{align*}
and $F\left( P(t),  {P}'(t)\right) = f\left( P(t) \right) + g\left( {P}'(t) \right)$.

Similar to the leapfrog scheme, our low-rank approach is designed based on an equivalent first-order system.
Introducing an auxiliary variable $Q(t) = P'(t)$, Eq.~\eqref{eq2.1} is transformed to 
the following first-order system:
\begin{equation}\label{eq2.2}
	\begin{split}
		\begin{bmatrix}
			P(t)\\
			Q(t)
		\end{bmatrix}'
		=
		& \begin{bmatrix}
			Q(t) \\
			-\alpha \left( A_x P(t) + P(t) A_y \right) - \beta \left( A_x Q(t) + Q(t) A_y \right)
			- \delta P(t) - \gamma Q(t)
		\end{bmatrix} + \\
		& \begin{bmatrix}
			\bm{0} \\
			F\left( P(t),  Q(t)\right)
		\end{bmatrix}, \qquad 
		\begin{bmatrix}
			P(0)\\
			Q(0)
		\end{bmatrix} = 
		\begin{bmatrix}
			P^0\\
			Q^0
		\end{bmatrix},
	\end{split}
\end{equation}
where $\bm{0}$ is a zero matrix of appropriate size.
To obtain a low-rank solution of Eq.~\cref{eq1.1}, 
the DLR approach \cite{koch2007dynamical,lubich2014projector} 
can be applied to Eq.~\cref{eq2.2}.
However, some numerical results in \cite[Section 7.2.1]{schrammer2022dynamical} demonstrated that 
the quality of its low-rank solutions deteriorates over time.
One possible reason should be that 
this direct application ignores some inherent structure of Eq.~\cref{eq1.1}.
Thus, we have to consider some other strategies to overcome this.

\section{The low-rank approximation}
\label{sec3}

In this section, we derive our DLR approximation to Eq.~\cref{eq1.1}
based on the first-order matrix differential system \cref{eq2.2}.

We first split Eq.~\cref{eq2.2} into three subproblems:
\begin{equation}\label{eq3.1}
	\begin{bmatrix}
		P_x(t)\\
		Q_x(t)
	\end{bmatrix}'
	=
	\begin{bmatrix}
		\bm{0} & \omega_1 I_x \\
		-\alpha A_x - \frac{1}{2} \delta I_x & -\beta A_x - \frac{1}{2} \gamma I_x
	\end{bmatrix}
	\begin{bmatrix}
		P_x(t)\\
		Q_x(t)
	\end{bmatrix}, \qquad
	\begin{bmatrix}
		P_x(0)\\
		Q_x(0)
	\end{bmatrix} = 
	\begin{bmatrix}
		P_x^0\\
		Q_x^0
	\end{bmatrix},
\end{equation}
\begin{equation}\label{eq3.2}
\begin{split}
	 \begin{bmatrix}
		P_y(t) & Q_y(t)
	\end{bmatrix}'& = 
	\begin{bmatrix}
		P_y(t) & Q_y(t)
	\end{bmatrix}
	\begin{bmatrix}
		\bm{0} & -\alpha A_y - \frac{1}{2} \delta I_y \\
		\omega_2 I_y & -\beta A_y - \frac{1}{2} \gamma I_y
	\end{bmatrix}, \\
    & \begin{bmatrix}
		P_y(0) & Q_y(0)
	\end{bmatrix} = 
	\begin{bmatrix}
		P_y^0 & Q_y^0
	\end{bmatrix}
\end{split}
\end{equation}
and 
\begin{equation}\label{eq3.3}
	\begin{bmatrix}
		P_{F}(t)\\
		Q_{F}(t)
	\end{bmatrix}' =
	\begin{bmatrix}
		\omega_3 Q_{F}(t)\\
		F\left( P_{F}(t),  Q_{F}(t) \right)
	\end{bmatrix}, \qquad
	\begin{bmatrix}
		P_{F}(0)\\
		Q_{F}(0)
	\end{bmatrix} = 
	\begin{bmatrix}
		P_{F}^0 \\
		Q_{F}^0
	\end{bmatrix},
\end{equation}
where $\omega_1 + \omega_2 + \omega_3 = 1$, $\omega_i > 0 (i = 1,2,3)$,
$I_x \in \mathbb{R}^{(N_x - 1) \times (N_x - 1)}$ and $I_y \in \mathbb{R}^{(N_y - 1) \times (N_y - 1)}$
are two identity matrices.
For a given positive integer $M$, let $\tau = \frac{T}{M}$ and $t_k = k \tau (k = 0,1,\cdots,M)$.
Denote by 
\begin{equation*}
	\Phi_{\tau}^{x} \left(\begin{bmatrix} P_x^0 \\ Q_x^0 \end{bmatrix} \right), \quad 
	\Phi_{\tau}^{y} \left(\begin{bmatrix} P_y^0 \\ Q_y^0 \end{bmatrix} \right)
	\quad \mathrm{and} \quad
	\Phi_{\tau}^{F} \left(\begin{bmatrix} P_{F}^0 \\ Q_{F}^0 \end{bmatrix} \right)
\end{equation*}
the solutions of the subproblems \cref{eq3.1}-\cref{eq3.3} with time step size $\tau$, respectively.
Applying a Strang splitting to Eq.~\cref{eq2.2}, we have the full-rank splitting scheme
\begin{equation} \label{eq3.4}
	\mathcal{L}_{\tau} = \Phi_{\tau/2}^{x} \circ \Phi_{\tau/2}^{y} \circ \Phi_{\tau}^{F} 
	\circ \Phi_{\tau/2}^{y} \circ \Phi_{\tau/2}^{x}.
\end{equation}
Then, starting with 
\begin{equation*}
	\begin{bmatrix} P_x^0 \\ Q_x^0 \end{bmatrix} = 
	\begin{bmatrix} P^0 \\ Q^0 \end{bmatrix},
\end{equation*}
the numerical solution of Eq.~\cref{eq2.2} at $t_k (k = 1,2,\cdots,M)$ is computed by
\begin{equation*}
	\begin{bmatrix}
		P^k \\
		Q^k
	\end{bmatrix} = 
	\mathcal{L}_{\tau}^k \left(\begin{bmatrix} P^0 \\ Q^0 \end{bmatrix} \right),
\end{equation*}
where $P^k$ and $Q^k$ are numerical approximations of $P(t_k)$ and $Q(t_k)$, respectively.

It is known that compared with standard numerical methods such as implicit Runge–Kutta methods, 
the DLR approximation \cite{koch2007dynamical} is an efficient numerical method 
for large-scale first-order matrix differential equations. 
Its key idea is to project the right-hand side of the first-order matrix differential equation 
onto the tangent space of the manifold of rank-$r$ matrices, where $r$ is small.
This yields an optimization problem. It can be solved by using
the projector-splitting integrator \cite{lubich2014projector}
or the unconventional robust integrator \cite{Ceruti2022unconventional}. 
Combining the splitting scheme Eq.~\cref{eq3.4} and the DLR approximation, 
we now show the details of finding the low-rank solution of Eq.~\cref{eq1.1}.

Let 
\begin{equation*}
	\mathcal{M}_{r} = \left\{ Y \in \mathbb{R}^{(N_x - 1) \times (N_y - 1)} 
	\mid \mathrm{rank} \left( Y \right) = r \ll \min \left\{ N_x - 1, N_y - 1 \right\} \right\}
\end{equation*}
be the manifold of rank-$r$ matrices and 
$
\mathcal{V}_{m,r} = \left\{ U \in \mathbb{R}^{m \times r}  
\mid U^{\top} U = I_r \right\}
$
be the Stiefel manifold.
Here $I_r$ is the identity matrix with order $r$. 
We know from \cite{lubich2014projector} that
for any matrix $Y \in \mathcal{M}_{r}$ can be expressed as
$Y = U S V^{\top}$,
where $U \in \mathcal{V}_{(N_x - 1),r}$, $V \in \mathcal{V}_{(N_y - 1),r}$ 
and $S$ is invertible and contains all singular values of $Y$.
It is worth mentioning that such expression is similar to the (truncated) singular value decomposition.
The distinction between them is that $S$ does not have to be a diagonal matrix in our case.

\subsection{Low-rank solutions to the two linear subproblems}
\label{sec3.1}

For the subproblem \cref{eq3.1}, 
we seek low-rank approaches $\tilde{P}_x(t) \in \mathcal{M}_{r_P}$ 
and $\tilde{Q}_x(t) \in \mathcal{M}_{r_Q}$ to $P_x(t)$ and $Q_x(t)$, respectively.
They satisfy the following equation:
\begin{equation} \label{eq3.5}
	\begin{bmatrix}
		\tilde{P}_x(t)\\
		\tilde{Q}_x(t)
	\end{bmatrix}'
	=
	\begin{bmatrix}
		\bm{0} & \omega_1 I_x \\
		-\alpha A_x - \frac{1}{2} \delta I_x & -\beta A_x - \frac{1}{2} \gamma I_x
	\end{bmatrix}
	\begin{bmatrix}
		\tilde{P}_x(t)\\
		\tilde{Q}_x(t)
	\end{bmatrix}, \qquad
	\begin{bmatrix}
		\tilde{P}_x(0)\\
		\tilde{Q}_x(0)
	\end{bmatrix} = 
	\begin{bmatrix}
		\tilde{P}_x^0\\
		\tilde{Q}_x^0
	\end{bmatrix},
\end{equation}
where $\tilde{P}_x^0 = U_0^{x} S_0^{x} (V_0^{x})^{\top} \in \mathcal{M}_{r_P}$ 
and $\tilde{Q}_x^0 = R_0^{x} \Sigma_0^{x} (W_0^{x})^{\top} \in \mathcal{M}_{r_Q}$ are 
rank-$r_P$ and rank-$r_Q$ approximations of $P_x^0$ and $Q_x^0$, respectively. 
The exact solution of this equation at $t_1$ can be expressed as
\begin{equation*}
	\begin{bmatrix}
		\tilde{P}_x(t_1) \\
		\tilde{Q}_x(t_1)
	\end{bmatrix}
	= B^x_{\tau}
	\begin{bmatrix}
		\tilde{P}_x^0\\
		\tilde{Q}_x^0
	\end{bmatrix},
\end{equation*}
where 
\begin{equation*}
	B^x_{\tau} = 
	\exp \left( \tau 
	\begin{bmatrix}
		\bm{0} & \omega_1 I_x \\
		-\alpha A_x - \frac{1}{2} \delta I_x & -\beta A_x - \frac{1}{2} \gamma I_x
	\end{bmatrix} \right).
\end{equation*}

Unfortunately, it cannot preserve the ranks of 
the input initial values $\tilde{P}_x^0$ and $\tilde{Q}_x^0$.
Thus, we apply the DLR approximation \cite{koch2007dynamical} to Eq.~\cref{eq3.5}.
We denote 
\begin{equation*}
	\begin{bmatrix} \tilde{P}_x^1 \\ \tilde{Q}_x^1 \end{bmatrix} = 
	\tilde{\Phi}_{\tau}^{x} \left(\begin{bmatrix} \tilde{P}_x^0 \\ \tilde{Q}_x^0 \end{bmatrix} \right)
\end{equation*} 
be the low-rank approach of $\left[ \begin{smallmatrix} P_x(t) \\ Q_x(t) \end{smallmatrix} \right]$ at $t_1$.

With the help of these notations, 
$\left[ \begin{smallmatrix} \tilde{P}_x^1 \\ \tilde{Q}_x^1 \end{smallmatrix} \right]$
is computed via \cref{alg1}, i.e., 
\begin{equation*}
	\mathtt{DLR\_L}(U_0^{x}, S_0^{x}, V_0^{x}, R_0^{x}, \Sigma_0^{x}, W_0^{x},\Delta_1,\Delta_2),
	\quad \mathrm{where} 
	\begin{bmatrix} \Delta_1 \\ \Delta_2 \end{bmatrix}
	= B^x_{\tau}
	\begin{bmatrix} \tilde{P}_x^0 \\ \tilde{Q}_x^0 \end{bmatrix}.
\end{equation*}
\begin{algorithm}[th]
	\caption{$\mathtt{DLR\_L}(\tilde{U}_0, \tilde{S}_0, \tilde{V}_0,
		\tilde{R}_0, \tilde{\Sigma}_0, \tilde{W}_0,\Delta_1,\Delta_2)$}
	\label{alg1}
	\begin{algorithmic}[1]
	\STATE {Compute QR-decomposition: $\tilde{U}_1 \hat{S} = \Delta_1 \tilde{V}_0$}
	\STATE {Compute QR-decomposition: $\tilde{V}_1 \tilde{S}_1^{\top} = \Delta_1^{\top} \tilde{U}_1$}
	\STATE {Compute QR-decomposition: $\tilde{R}_1 \hat{S} = \Delta_2 \tilde{W}_0$}
	\STATE {Compute QR-decomposition: $\tilde{W}_1 \tilde{\Sigma}_1^{\top} = \Delta_2^{\top} \tilde{R}_1$}
	\STATE {Return: $\tilde{P}_1 = \tilde{U}_1 \tilde{S}_1 \tilde{V}_1^{\top}$ and 
		$\tilde{Q}_1 = \tilde{R}_1 \tilde{\Sigma}_1 \tilde{W}_1^{\top}$}
\end{algorithmic}
\end{algorithm}

Next, we turn to the subproblem \cref{eq3.2}.
Denote by $\tilde{P}_y^0 = U_0^{y} S_0^{y} (V_0^{y})^{\top} \in \mathcal{M}_{r_P}$ 
and $\tilde{Q}_y^0 = R_0^{y} \Sigma_0^{y} (W_0^{y})^{\top} \in \mathcal{M}_{r_Q}$ the 
rank-$r_P$ and rank-$r_Q$ approximations to the initial values $P_y^0$ and $Q_y^0$, respectively.
Similar to the subproblem \cref{eq3.1}, 
for given low-rank initial values $\tilde{P}_y^0$ and $\tilde{Q}_y^0$, 
low-rank solutions of the subproblem \cref{eq3.2} are obtained by solving the following equation: 
\begin{align*}
	\begin{bmatrix}
		\tilde{P}_y(t) & \tilde{Q}_y(t)
	\end{bmatrix}' &=
	\begin{bmatrix}
		\tilde{P}_y(t) & \tilde{Q}_y(t)
	\end{bmatrix}
	\begin{bmatrix}
		\bm{0} & -\alpha A_y - \frac{1}{2} \delta I_y \\
		\omega_2 I_y & -\beta A_y - \frac{1}{2} \gamma I_y
	\end{bmatrix}, \\
&	\begin{bmatrix}
		\tilde{P}_y(0) & \tilde{Q}_y(0)
	\end{bmatrix} = 
	\begin{bmatrix}
		\tilde{P}_y^0 & \tilde{Q}_y^0
	\end{bmatrix},
\end{align*}
where $\tilde{P}_y(t) \in \mathcal{M}_{r_P}$ 
and $\tilde{Q}_y(t) \in \mathcal{M}_{r_Q}$ are low-rank approximations to $P_x(t)$ and $Q_x(t)$, respectively.

Let  
\begin{equation*}
	B^y_{\tau} = 
	\exp \left( \tau 
	\begin{bmatrix}
		\bm{0} & -\alpha A_y - \frac{1}{2} \delta I_y \\
		\omega_2 I_y & -\beta A_y - \frac{1}{2} \gamma I_y
	\end{bmatrix} \right).
\end{equation*}
Then, the exact solution of this matrix differential equation at $t_1$ is 
\begin{equation*}
	\begin{bmatrix}
		\tilde{P}_y(t_1) & \tilde{Q}_y(t_1)
	\end{bmatrix}
	= \begin{bmatrix}
		\tilde{P}_y^0 & \tilde{Q}_y^0
	\end{bmatrix} B^y_{\tau}.
\end{equation*}
This is also not rank-preserving. Thus, we compute
the low-rank approximation 
\begin{equation*}
	\begin{bmatrix} \tilde{P}_y^1 & \tilde{Q}_y^1 \end{bmatrix} = 
	\tilde{\Phi}_{\tau}^{y} \left(\begin{bmatrix} \tilde{P}_y^0 & \tilde{Q}_y^0 \end{bmatrix} \right)
\end{equation*} 
of $\begin{bmatrix} P_y(t_1) & Q_y(t_1) \end{bmatrix}$ via \cref{alg1}, i.e., 
$\mathtt{DLR\_L}(U_0^{y}, S_0^{y}, V_0^{y}, R_0^{y}, \Sigma_0^{y}, W_0^{y},\Delta_1,\Delta_2)$,
where 
\begin{equation*}
	\begin{bmatrix} \Delta_1 & \Delta_2 \end{bmatrix}
	= \begin{bmatrix} \tilde{P}_y^0 &\ \tilde{Q}_y^0 \end{bmatrix} B^y_{\tau}.
\end{equation*}

Now, the remaining part of our low-rank algorithm is 
to find a low-rank solution to the nonlinear subproblem \cref{eq3.3}.
This is more difficult than subproblems \cref{eq3.1} and \cref{eq3.2}.
In the next subsection, we provide a method for finding the low-rank solution of Eq.~\cref{eq3.3}.

\subsection{A low-rank solution to the nonlinear subproblem}
\label{sec3.2}

Compared to subproblems \cref{eq3.1} and \cref{eq3.2}, 
finding a low-rank solution for the subproblem \cref{eq3.3} is more complicated.
Firstly, we split the nonlinear part \cref{eq3.3} into following two subproblems:
\begin{equation*}
	\begin{bmatrix}	P_{FL}(t) \\ Q_{FL}(t) \end{bmatrix}' =
	\begin{bmatrix}	\omega_3 Q_{FL}(t) \\ \bm{0} \end{bmatrix}, \qquad
	\begin{bmatrix}	P_{FL}(0) \\ Q_{FL}(0) \end{bmatrix} = 
	\begin{bmatrix} P_{FL}^0 \\ Q_{FL}^0 \end{bmatrix},
\end{equation*}
\begin{equation*}
	\begin{bmatrix} P_{FN}(t) \\ Q_{FN}(t) \end{bmatrix}' =
	\begin{bmatrix} \bm{0} \\ F\left( P_{FN}(t),  Q_{FN}(t) \right) \end{bmatrix}, \qquad
	\begin{bmatrix} P_{FN}(0) \\ Q_{FN}(0) \end{bmatrix} = 
	\begin{bmatrix} P_{FN}^0 \\ Q_{FN}^0 \end{bmatrix}.
\end{equation*}

Apply a Strang splitting combined with the DLR approximation to these subproblems.
Then, we obtain the low-rank solutions $\tilde{P}_F^1$ and $\tilde{Q}_F^1$ of \cref{eq3.3} at $t_1$ by
starting from low-rank initial values 
\begin{equation*}
	\tilde{P}_F^0 = U_0^{F} S_0^{F} (V_0^{F})^{\top} \in \mathcal{M}_{r_P} \approx P_F^0 = P_{FN}^0
\end{equation*}
and
\begin{equation*}
	\tilde{Q}_F^0 = R_0^{F} \Sigma_0^{F} (W_0^{F})^{\top} \in \mathcal{M}_{r_Q} \approx Q_F^0 = Q_{FN}^0,
\end{equation*} 
see \cref{alg2}. Moreover, we denote 
$\left[
\begin{smallmatrix} \tilde{P}_F^1 \\ \tilde{Q}_F^1 \end{smallmatrix} \right] = 
\tilde{\Phi}_{\tau}^{F} \left(\left[ \begin{smallmatrix} \tilde{P}_F^0 \\ \tilde{Q}_F^0 \end{smallmatrix} \right] \right).
$
\begin{algorithm}[ht] 
	\caption{Calculate low-rank approximations to $P_F(t)$ and $Q_F(t)$, single time step}
	\label{alg2}
	\begin{algorithmic}[1]
		\STATE {Input: $U_0^{F}$, $S_0^{F}$, $V_0^{F}$, $R_0^{F}$, $\Sigma_0^{F}$, $W_0^{F}$
			and $\omega_3$}
		\STATE {Compute $R^{F,1}$, $\Sigma^{F,1}$ and $W^{F,1}$
			for $t \in[0,\frac{\tau}{2}]$ via \cref{alg3}, i.e.,  
			\begin{equation*}
				\mathtt{DLR\_FN}(U_0^{F},S_0^{F},V_0^{F},R_0^{F},\Sigma_0^{F},W_0^{F},F,\tau/2)
		\end{equation*}}\vspace{-4mm}
		\STATE {Compute $U^{F,1}$, $S^{F,1}$ and $V^{F,1}$
			for $t \in[0,\tau]$ via \cref{alg4}, i.e.,
			\begin{equation*}
				\mathtt{DLR\_FL}(U_0^{F},S_0^{F},V_0^{F},R^{F,1},\Sigma^{F,1},W^{F,1},\tau,\omega_3)
		\end{equation*}}\vspace{-4mm}
		\STATE {Compute $R^{F,2}$, $\Sigma^{F,2}$ and $W^{F,2}$
			for $t \in[0,\frac{\tau}{2}]$ via \cref{alg3}, i.e.,  
			\begin{equation*}
				\mathtt{DLR\_FN}(U^{F,1},S^{F,1},V^{F,1},R^{F,1},\Sigma^{F,1},W^{F,1},F,\tau/2)
		\end{equation*}}\vspace{-4mm}
		\STATE {Output: $\tilde{P}_F^1 = U^{F,1} S^{F,1} (V^{F,1})^{\top}$ and 
			$\tilde{Q}_F^1 = R^{F,2} \Sigma^{F,2} (W^{F,2})^{\top}$}
	\end{algorithmic}
\end{algorithm}
\begin{algorithm}[H]
	\caption{$\mathtt{DLR\_FN}(\tilde{U}_0, \tilde{S}_0, \tilde{V}_0, \tilde{R}_0, 
		\tilde{\Sigma}_0, \tilde{W}_0, F, \tau)$}
\label{alg3}
\begin{algorithmic}[1]
\IF {$g(u_t) \neq 0$}
\STATE {Compute $K_1 = K(\tau)$ by solving
	\begin{equation*}
		K'(t) = F \left( \tilde{U}_0 \tilde{S}_0 \tilde{V}_0^\top, K(t) \tilde{W}_0^\top \right) \tilde{W}_0, \qquad 
		K(0) = \tilde{R}_0 \tilde{\Sigma}_0,~t \in [0,\tau]
\end{equation*}}\vspace{-4mm}
\STATE {Compute QR-decomposition: $\tilde{R}_1 \hat{\Sigma}_0 = K_1$}
\STATE {Compute $\hat{\Sigma}_1 = \hat{\Sigma}(\tau)$ by solving
	\begin{align*}
		\hat{\Sigma}'(t) = F \left( \tilde{U}_0 \tilde{S}_0 \tilde{V}_0^\top, 
		\tilde{R}_1 \hat{\Sigma}(t) \tilde{W}_0^\top \right) \tilde{W}_0, \qquad 
		\hat{\Sigma}(0) = \hat{\Sigma}_0,~t \in [0,\tau]
\end{align*}}\vspace{-4mm}
\STATE {Compute $L_1 = L(\tau)$ by solving
	\begin{align*}
		L'(t) = F \left( \tilde{U}_0 \tilde{S}_0 \tilde{V}_0^\top, 
		\tilde{R}_1 L^\top(t) \right)^\top \tilde{R}_1, \qquad 
		L(0) = \tilde{W}_0 \hat{\Sigma}_1^\top,~t \in [0,\tau]
\end{align*}}\vspace{-4mm}
\STATE {Compute QR-decomposition: $\tilde{W}_1 \tilde{\Sigma}_1^{\top} = L_1$}
\ELSE 
\STATE {\% $g(u_t) = 0$}
\STATE {$\Delta Q = \tau F \left(\tilde{U}_0 \tilde{S}_0 \tilde{V}_0^\top \right)$,
	$\Delta K = \Delta Q\, \tilde{W}_0$}
\STATE {$K_1 = \tilde{R}_0 \tilde{\Sigma}_0 + \Delta K$}
\STATE {Compute QR-decomposition: $\tilde{R}_1 \hat{\Sigma}_0 = K_1$}
\STATE {$\hat{\Sigma}_1 = \hat{\Sigma}_0 - \tilde{R}_1^\top \Delta K$}
\STATE {$L_1 = \tilde{W}_0 \hat{\Sigma}_1^\top + \Delta Q^\top \tilde{R}_1$}
\STATE {Compute QR-decomposition: $\tilde{W}_1 \tilde{\Sigma}_1^{\top} = L_1$}
\ENDIF
\STATE {Return: $\tilde{R}_1$, $\tilde{\Sigma}_1$ and $\tilde{W}_1$}
\end{algorithmic}
\end{algorithm}

Let $\tilde{P}^k = U_k S_k V_k^\top \in \mathcal{M}_{r_P}$ and 
$\tilde{Q}^k = R_k \Sigma_k W_k^\top \in \mathcal{M}_{r_Q}$ be rank-$r_P$ 
and rank-$r_Q$ approaches of $P(t_k)$ and $Q(t_k)$ ($k = 0,1,\cdots,M$), respectively.
Starting from  $\tilde{P}_x^0 = \tilde{P}^0$ and $\tilde{Q}_x^0 = \tilde{Q}^0$, 
we obtain a low-rank approach to Eq.~\cref{eq2.2} at $t_k$, i.e., 
\begin{equation*}
	\begin{bmatrix}
		\tilde{P}^k \\
		\tilde{Q}^k
	\end{bmatrix} = 
	\mathcal{L}_{\tau,r}^k \left(\begin{bmatrix} \tilde{P}^0 \\ \tilde{Q}^0 \end{bmatrix} \right), \qquad
	k = 1, \cdots,M,
\end{equation*}
where 
\begin{equation}\label{eq3.6}
	\mathcal{L}_{\tau,r} = \tilde{\Phi}_{\tau/2}^{x} \circ \tilde{\Phi}_{\tau/2}^{y} \circ \tilde{\Phi}_{\tau}^{F} 
	\circ \tilde{\Phi}_{\tau/2}^{y} \circ \tilde{\Phi}_{\tau/2}^{x}
\end{equation}
is our low-rank scheme.

In the next section, we provide several numerical examples to 
test the low-rank scheme \cref{eq3.6}.
\begin{algorithm}[ht]
	\caption{$\mathtt{DLR\_FL}(\tilde{U}_0, \tilde{S}_0, \tilde{V}_0, \tilde{R}_0, 
		\tilde{\Sigma}_0, \tilde{W}_0,\tau,\omega)$}
\label{alg4}
\begin{algorithmic}[1] 
	\STATE {$\Delta P = \tau \omega \tilde{R}_0 \tilde{\Sigma}_0 \tilde{W}_0^\top$, 
		$\Delta K = \Delta P\, \tilde{V}_0$}	
	\STATE {$K_1 = \tilde{U}_0 \tilde{S}_0 + \Delta K$}
	\STATE {Compute QR-decomposition: $\tilde{U}_1 \hat{S}_0 = K_1$}
	\STATE {$\hat{S}_1 = \hat{S}_0 - \tilde{U}_1^\top \Delta K$}
	\STATE {$L_1 = \tilde{V}_0 \hat{S}_1^\top + \Delta P^\top \tilde{U}_1$}
	\STATE {Compute QR-decomposition: $\tilde{V}_1 \tilde{S}_1^{\top} = L_1$}
	\STATE {Return: $\tilde{U}_1$, $\tilde{S}_1$ and $\tilde{V}_1$}
\end{algorithmic}
\end{algorithm}

\section{Numerical experiments}
\label{sec4}

In this section, we report two examples to show the convergence order of \cref{eq3.6}.
Example 3 simulates the time evolution of Eq.~\cref{eq1.1} with different initial values
to explore the performance of our scheme \cref{eq3.6}.
Denote
\begin{equation*}
	\textrm{relerr}(\tau) = \frac{\left\| \tilde{P}^M - \mathcal{U}(T) \right\|_{F}}{\left\| \mathcal{U}(T) \right\|_{F}} \quad \mathrm{and} \quad
	\textrm{rate} = \log_{\tau_1/\tau_2} \frac{\textrm{relerr}(\tau_1)}{\textrm{relerr}(\tau_2)}, 
\end{equation*}
where $\mathcal{U}(T) = \left[ u(x_i,y_j,T) \right]_{\substack{1 \leq i \leq N_x - 1 \\ 
		1 \leq j \leq N_y - 1}}$
and $\| \cdot \|_{F}$ is the Frobenius norm.
For clarity, in this work, we fix $N_x = N_y = N$ and $r_P = r_Q = r$.

All experiments were done by MATLAB R2018b on a Windows 10 (64 bit) i5-1135G7 CPU 2.40GHz, 
24 GB of RAM.

\noindent\textbf{Example 1.} In Eq.~\cref{eq1.1}, 
we consider $\alpha = 1$, $\beta = 0.1$, $\gamma = 0.001$,
$\delta = 1$, $\Omega = [0,1]^2$, the initial values
\begin{equation*}
	p(x,y) = 2 \sin(3 \pi x) \sin(3 \pi y),\qquad  q(x,y) = -\sin(3 \pi x) \sin(3 \pi y)
\end{equation*}
and the nonlinear terms
\begin{equation*}
	f(u) = u^2,\qquad g(v) = \sin(v).
\end{equation*}

In this example, we use the numerical solution ($(N,M) = (512,80000)$) computed 
by the EI-SW4 method \cite{phan2022exponential} as the reference solution
since the exact solution is unknown.

\begin{figure}[ht]
	\centering
	{\includegraphics[width=2.5in,height=2.5in]{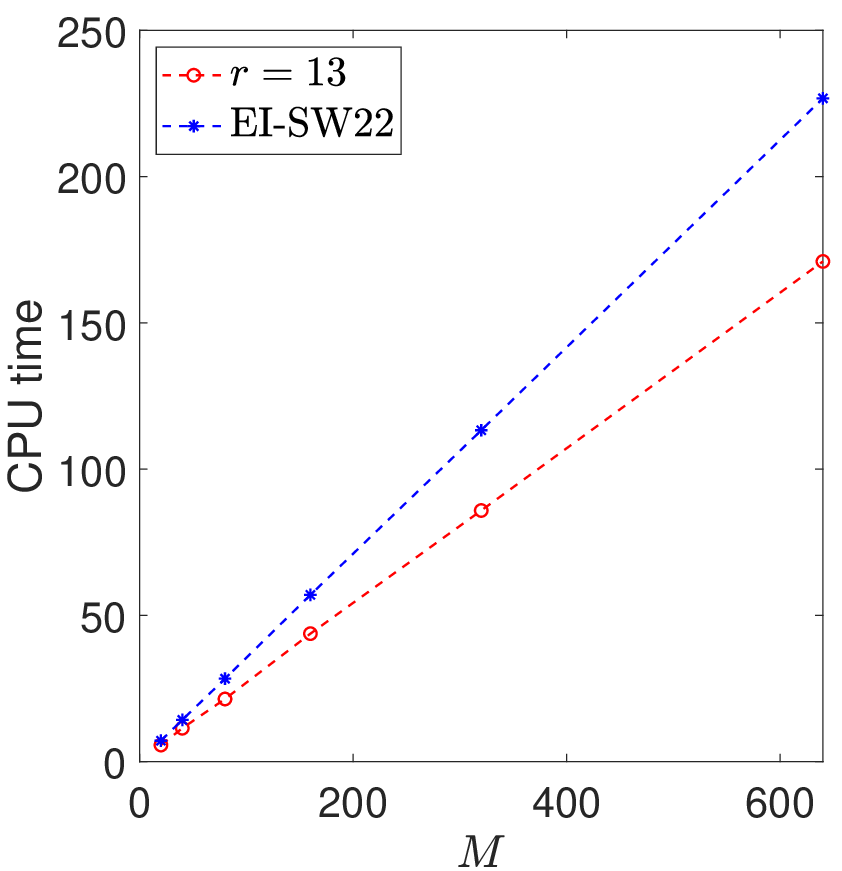}}
	\caption{Comparison of the CPU time (in seconds) of our method and the EI-SW22 method \cite{phan2022exponential} for Example 1, 
		where $N = 512$ and $M \in \{20,40, \cdots, 20 \cdot 2^{5}\}$.}
	\label{fig1}
\end{figure}

We report the relative errors and time convergence order of the scheme \cref{eq3.6} in \cref{tab1}.
From this table, for $r = 11,13 (\ll N)$, the time convergence order is $2$ as expected.
For $r = 7,9$ and large $M$, the errors start to stagnate.
One possible explanation for this phenomenon is 
that the rank truncation error dominates the global error.
\Cref{fig1} compares the CPU time required 
by the proposed scheme \cref{eq3.6} and the EI-SW22 method \cite{phan2022exponential}
for $M \in \{20,40, \cdots, 20 \cdot 2^{5}\}$.
Clearly, our method needs less CPU time to solve Eq.~\cref{eq1.1}.
\begin{table}[ht]\tabcolsep=2.0pt
	\caption{Errors and time convergence orders for Example 1, where $T = 0.1$, $N = 512$ 
		and $(\omega_1, \omega_2, \omega_3) = (0.98, 0.01, 0.01)$ .}
	\centering
	\begin{tabular}{ccccccccc}
		\hline
		& \multicolumn{2}{c}{$r = 7$} & \multicolumn{2}{c}{$r = 9$} & \multicolumn{2}{c}{$r = 11$} 
		& \multicolumn{2}{c}{$r = 13$} \\
		[-2pt] \cmidrule(lr){2-3} \cmidrule(lr){4-5} \cmidrule(lr){6-7} \cmidrule(lr){8-9} 
		\\ [-11pt]
		$M$ & $\textrm{relerr}(\tau)$ & $\textrm{rate}$ 
		& $\textrm{relerr}(\tau)$ & $\textrm{rate}$
		& $\textrm{relerr}(\tau)$ & $\textrm{rate}$ & $\textrm{relerr}(\tau)$ & $\textrm{rate}$ \\
		\hline
		20 & 8.4712E-05& -- &		8.4749E-05&	--&	8.4760E-05& --&		8.4760E-05	& -- \\
		40&	2.1922E-05&	1.9502& 	2.1848E-05&	1.9557 &	2.1765E-05&	1.9614 &	2.1772E-05&	1.9609 \\
		80&	6.2412E-06&	1.8125 &	5.8979E-06&	1.8892 &	5.5095E-06&	1.9820 &	5.5198E-06&	1.9798 \\
		160& 3.2702E-06& 0.9324& 	2.5368E-06&	1.2172 &	1.3871E-06&	1.9898 &	1.3867E-06&	1.9930 \\
		320& 2.9992E-06& 0.1248& 	2.1609E-06&	0.2314 &	3.7595E-07&	1.8835 &	3.4784E-07&	1.9952 \\
		\hline
	\end{tabular}
	\label{tab1}
\end{table}

\noindent\textbf{Example 2.} Consider Eq.~\cref{eq1.1} in $\Omega = [0,1]^2$ with 
$\alpha = 1$, $\beta = 0.001$, $\gamma = 10^{-6}$,
$\delta = 1$, the initial values
\begin{equation*}
	p(x,y) = 10 \sin(3 \pi x) \sin(3 \pi y),\qquad  q(x,y) = -10 \cos(3 \pi x) \cos(3 \pi y)
\end{equation*}
and the nonlinear terms
\begin{equation*}
	f(u) = u^3,\qquad g(v) = 0.
\end{equation*}

The exact solution in this case is unknown. 
We choose the EI-SW4 method \cite{phan2022exponential} to compute 
the reference solution ($(N,M) = (512,80000)$).

\begin{figure}[ht]
	\centering
	{\includegraphics[width=2.5in,height=2.5in]{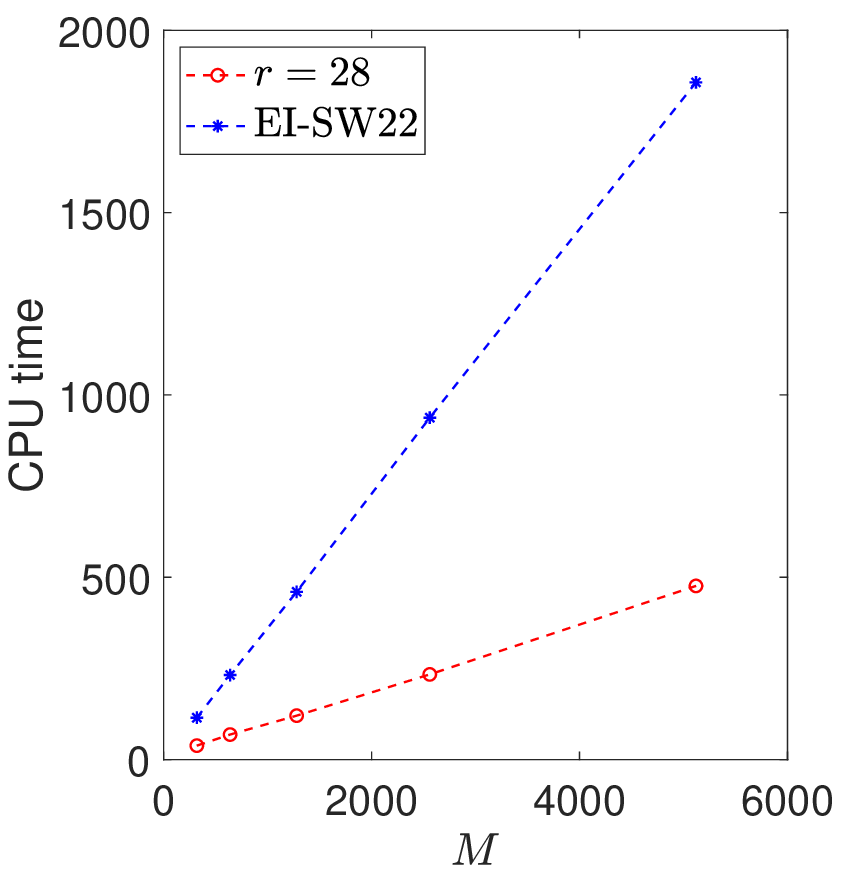}}
	\caption{Comparison of the CPU time (in seconds) of our method and the EI-SW22 method \cite{phan2022exponential} for Example 2, 
		where $N = 512$ and $M \in \{320,640, \cdots, 320 \cdot 2^{4}\}$.}
	\label{fig2}
\end{figure}

\Cref{tab2} lists the relative errors and observed convergence order in time.
We can see from this table that for $r = 28$, the time convergence order of our method is $2$.
For $r = 22,24,26$ and small $M = 320,640,1280$, the errors dropped as we expected.
However, for larger values of $M$, a stagnation of the error is observed. 
The reason may be that the rank truncation error 
prevails over the other components of the error.
\Cref{fig2} plots the CPU time of the proposed method and the EI-SW22 method \cite{phan2022exponential}
for $N = 512$ and $M \in \{320,640, \cdots, 320 \cdot 2^{4}\}$.
From this figure, the CPU time of the low-rank scheme \cref{eq3.6} is less than the EI-SW22 method.
To sum up, our low-rank scheme \cref{eq3.6} is efficient and robust for solving Eq.~\cref{eq1.1}.
\begin{table}[ht]\tabcolsep=2.0pt
	\caption{Errors and time convergence orders for $T = 1$, $N = 512$ 
		and $\omega_1=\omega_2=\omega_3 = \frac{1}{3}$ for Example 2.}
	\centering
	\begin{tabular}{ccccccccc}
		\hline
		& \multicolumn{2}{c}{$r = 22$} & \multicolumn{2}{c}{$r = 24$} & \multicolumn{2}{c}{$r = 26$} 
		& \multicolumn{2}{c}{$r = 28$} \\
		[-2pt] \cmidrule(lr){2-3} \cmidrule(lr){4-5} \cmidrule(lr){6-7} \cmidrule(lr){8-9} 
		\\ [-11pt]
		$M$ & $\textrm{relerr}(\tau)$ & $\textrm{rate}$ 
		& $\textrm{relerr}(\tau)$ & $\textrm{rate}$
		& $\textrm{relerr}(\tau)$ & $\textrm{rate}$ & $\textrm{relerr}(\tau)$ & $\textrm{rate}$ \\
		\hline
		320&	5.1735E-04&	--&	5.2065E-04&	--&	5.1557E-04&	--&	5.1553E-04&	-- \\
		640&	1.3561E-04&	1.9317& 1.3422E-04&	1.9557&	1.2897E-04&	1.9991& 1.2891E-04&	1.9997 \\
		1280&	4.4166E-05&	1.6185& 3.2789E-05&	2.0333&	3.2458E-05&	1.9904&	3.2274E-05&	1.9979 \\
		2560&	3.6424E-05&	0.2780& 1.1178E-05&	1.5525&	8.7559E-06&	1.8902&	8.1511E-06&	1.9853 \\
		5120&	2.3599E-05&	0.6262& 6.1876E-06&	0.8532&	4.5948E-06&	0.9303&	2.0215E-06&	2.0116 \\
		\hline
	\end{tabular}
	\label{tab2}
\end{table}

\noindent\textbf{Example 3.} In this example, we choose 
$\alpha = 0.6$, $\beta = 0.3$, $\gamma = 0.05$, $\delta = 0$ , $g(v) = 0$ and 
consider the following three different scenarios: \\
\textbf{5-petals flower}: $\Omega = [-3,3]^2$, $f(u) = u^2$, 
$(\omega_1, \omega_2, \omega_3) = (0.98,0.01,0.01)$ (for the scheme \cref{eq3.6}),
\begin{equation*}
	p(x,y) = 
	\begin{cases}
		0.1 \left(x^2 + y^2 + 1\right), 
		& x^2 + y^2 \leq \left[ \sin\left(5 \arctan \frac{y}{x}\right) + 1.5 \right]^2, \\
		0, & \mathrm{Others}
	\end{cases}
\end{equation*}
and $q(x,y) = 0.5\, p(x,y)$. \\
\textbf{Cardioid}: $\Omega = [-2.5,0.5] \times [-1.5,1.5]$, $f(u) = u (1 - u)$,
$\omega_1 = \omega_2 = \omega_3 = \frac{1}{3}$ (for the scheme \cref{eq3.6}),
\begin{equation*}
	p(x,y) = 
	\begin{cases}
		0.15 \exp \left( - \left(x^2 + y^2 + x\right)^2 + x^2 + y^2\right),
		& x^2 + y^2 + x \leq \sqrt{x^2 + y^2}, \\
		0, & \mathrm{Others}
	\end{cases}
\end{equation*}
and $ q(x,y) = -0.25\, p(x,y)$. \\
\textbf{Astroid}: $\Omega = [-1,1]^2$, $f(u) = |\sin(u)|$,
$(\omega_1, \omega_2, \omega_3) = (0.98,0.01,0.01)$ (for the scheme \cref{eq3.6}),
\begin{equation*}
	p(x,y) = 
	\begin{cases}
		- \left( x^{\frac{2}{3}} + y^{\frac{2}{3}} + 0.1 \right),
		& x^{\frac{2}{3}} + y^{\frac{2}{3}} \leq 0.7^{\frac{2}{3}}, \\
		0, & \mathrm{Others}
	\end{cases}
\end{equation*}
and $q(x,y) = 10\, p(x,y)$.

\begin{figure}[p]
	\setlength{\tabcolsep}{0.2pt}
	\centering
	\begin{tabular}{m{0.4cm}<{\centering} m{3.1cm}<{\centering} m{3.1cm}<{\centering} m{3.1cm}<{\centering} m{3.1cm}<{\centering}}
		& $t = 0$ & $t = 1$ & $t = 2$ & $t = 3$ \\
		\rotatebox{90}{EI-SW22} &
		\includegraphics[width=1.2in,height=1.2in]{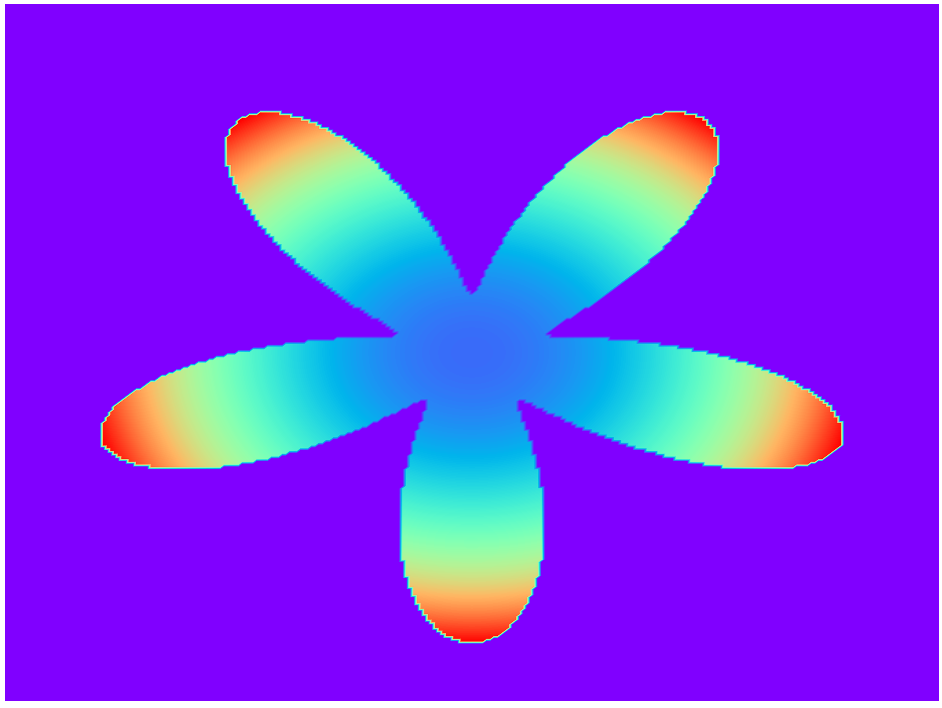} &
		\includegraphics[width=1.2in,height=1.2in]{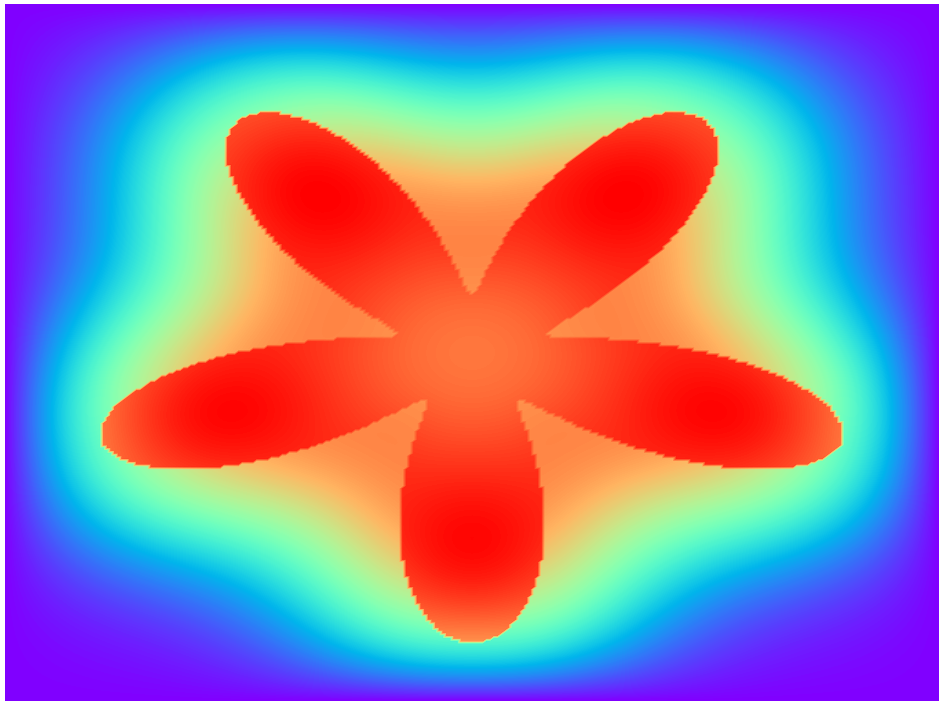} &
		\includegraphics[width=1.2in,height=1.2in]{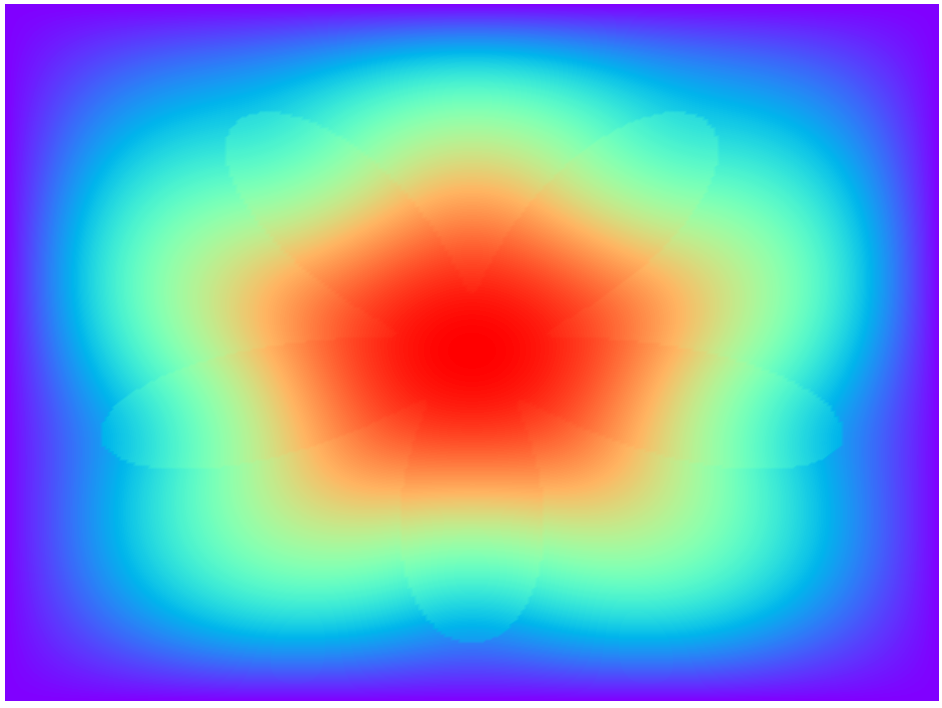} &
		\includegraphics[width=1.2in,height=1.2in]{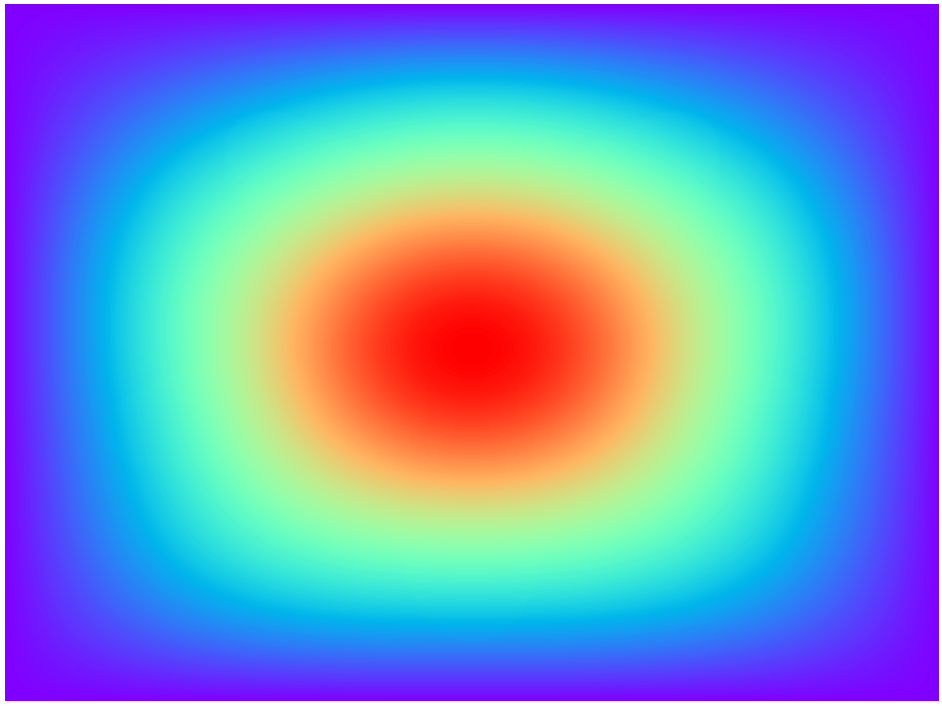} \\
		\rotatebox{90}{$r = 88$} &
		\includegraphics[width=1.2in,height=1.2in]{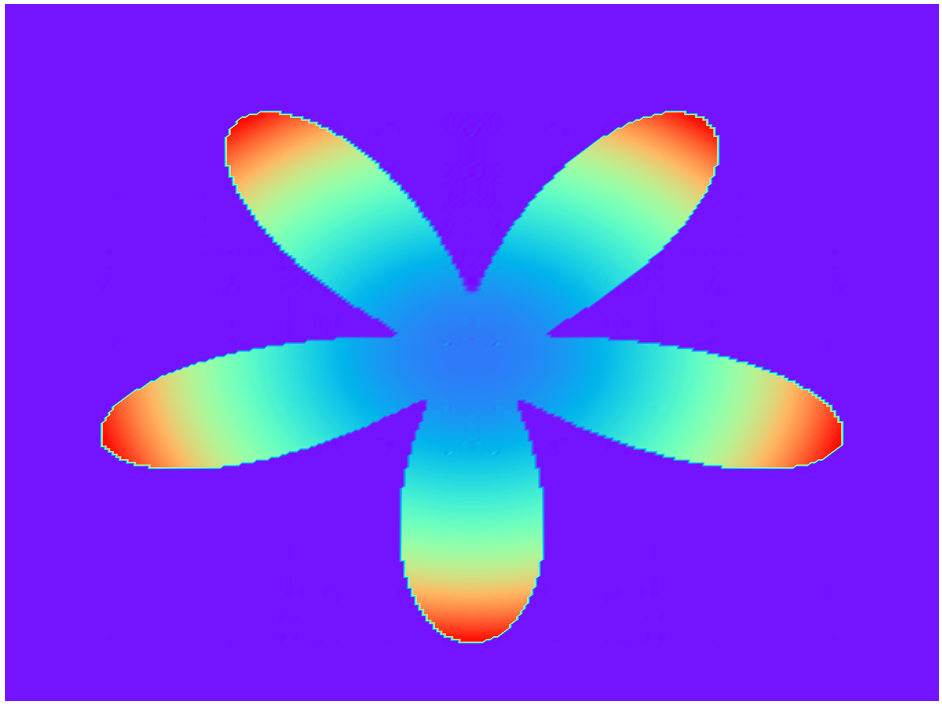} &
		\includegraphics[width=1.2in,height=1.2in]{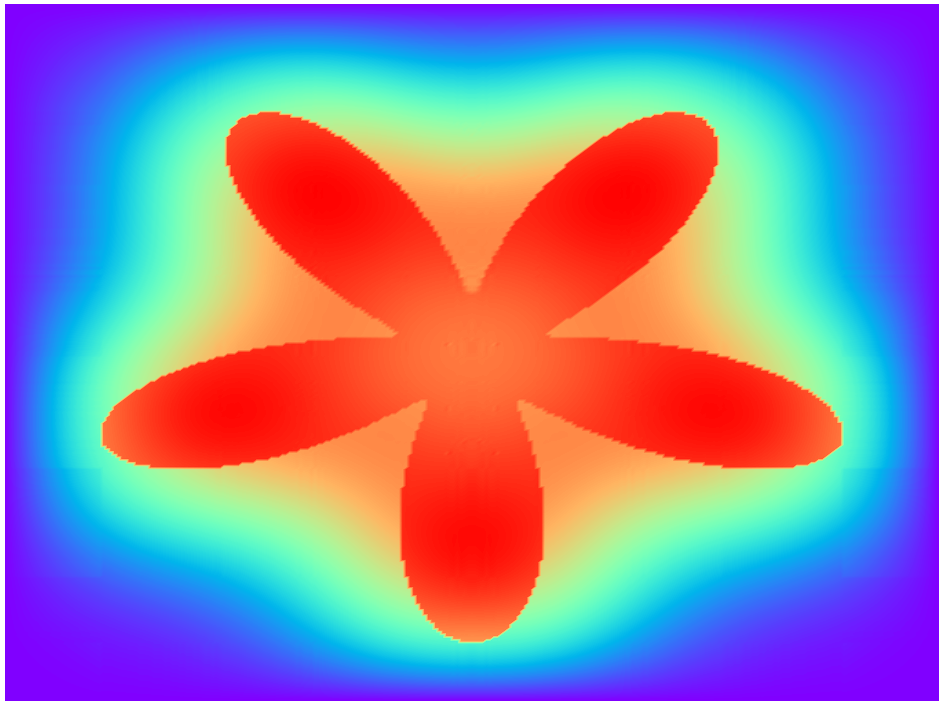} &
		\includegraphics[width=1.2in,height=1.2in]{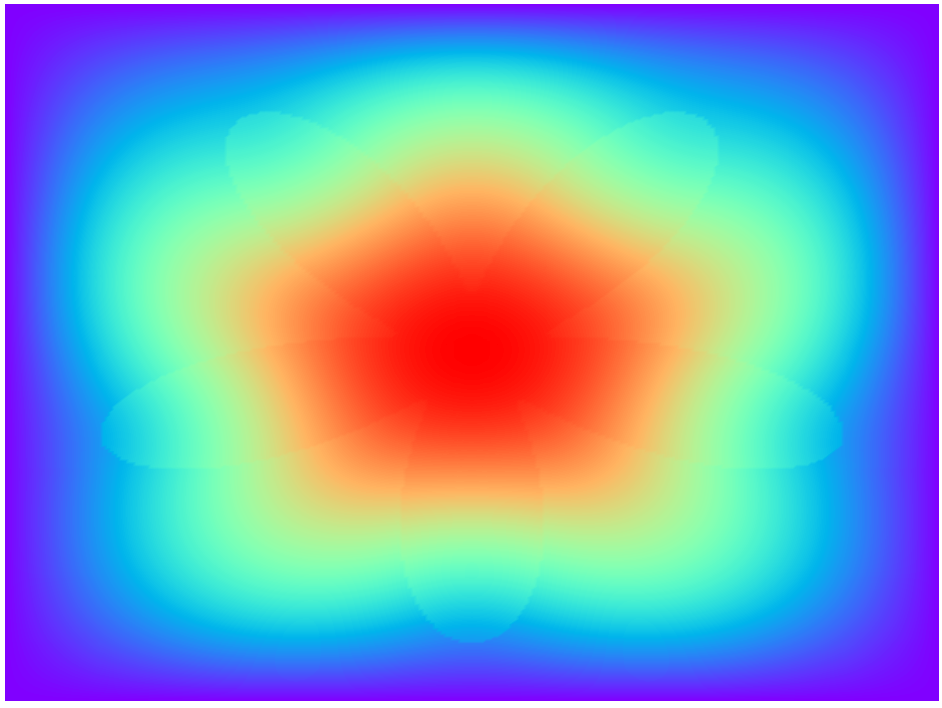} &
		\includegraphics[width=1.2in,height=1.2in]{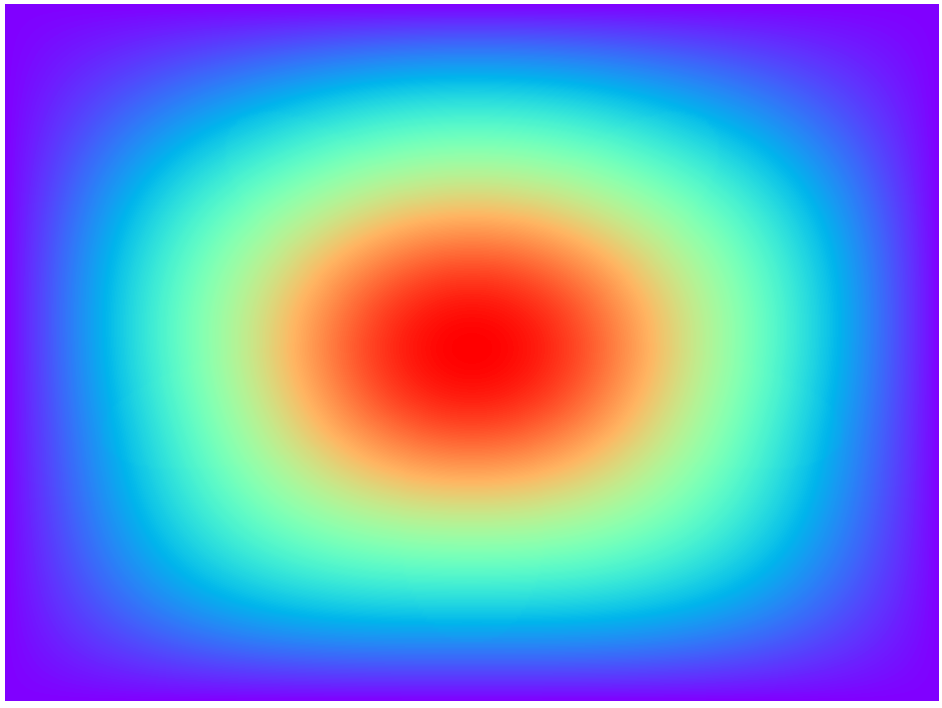} \\
		\rotatebox{90}{EI-SW22} &
		\includegraphics[width=1.2in,height=1.2in]{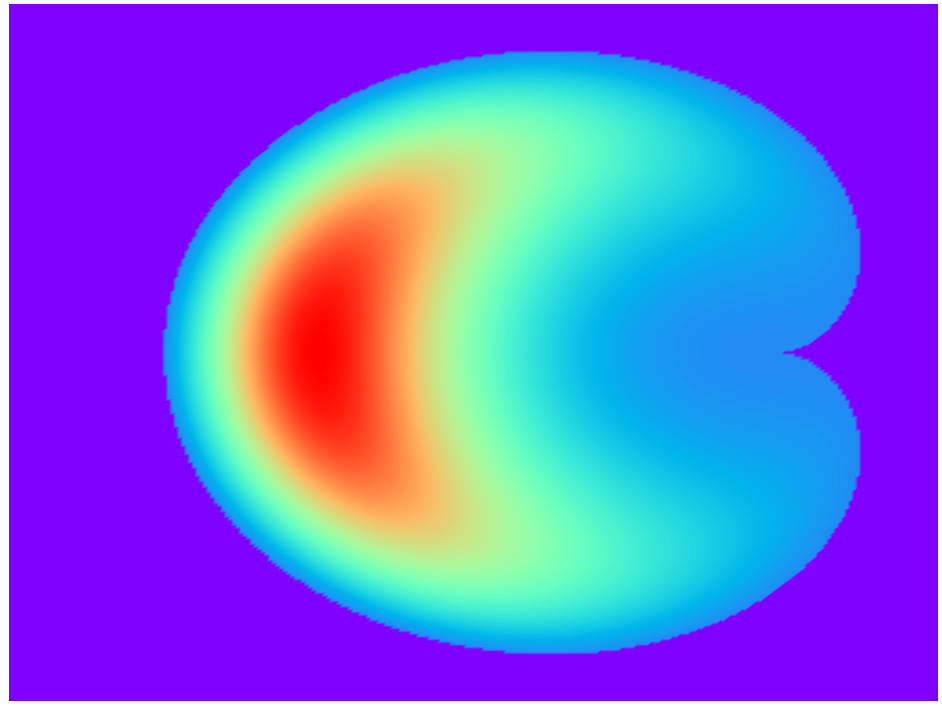} &
		\includegraphics[width=1.2in,height=1.2in]{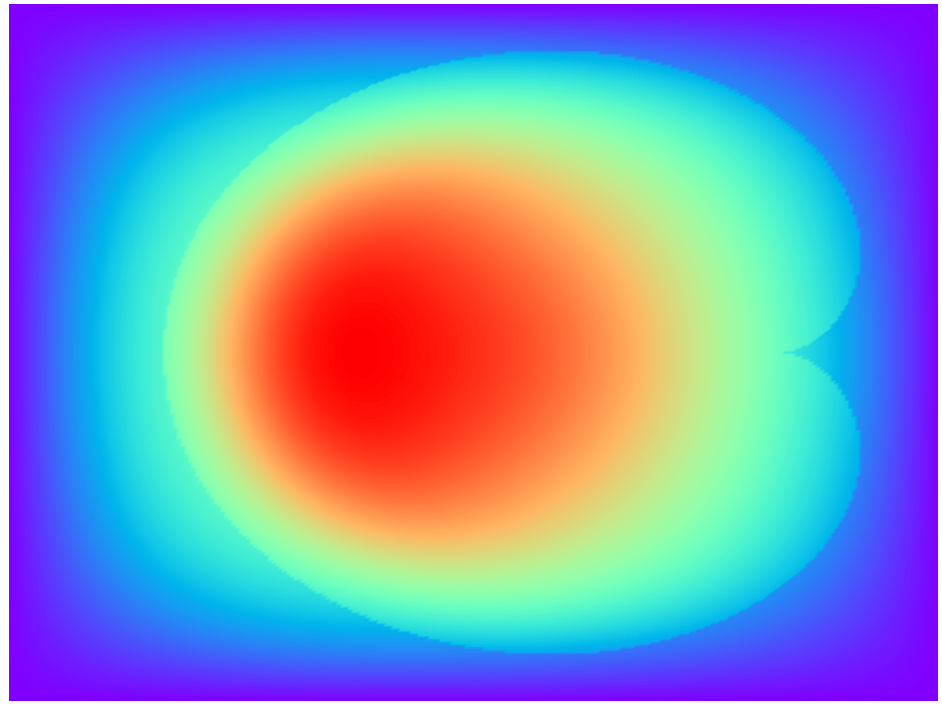} &
		\includegraphics[width=1.2in,height=1.2in]{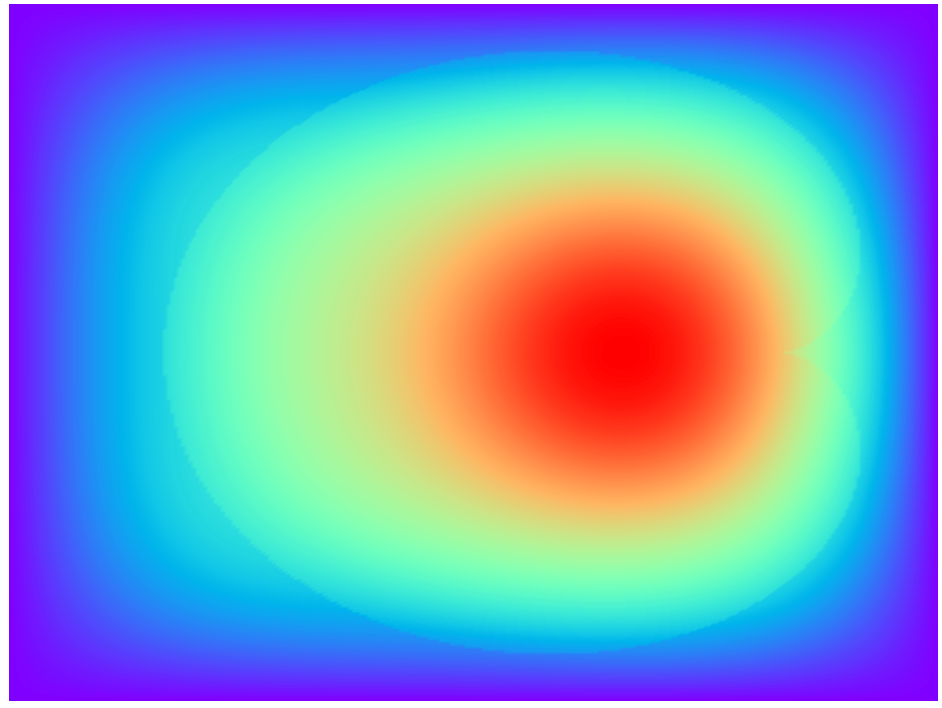} &
		\includegraphics[width=1.2in,height=1.2in]{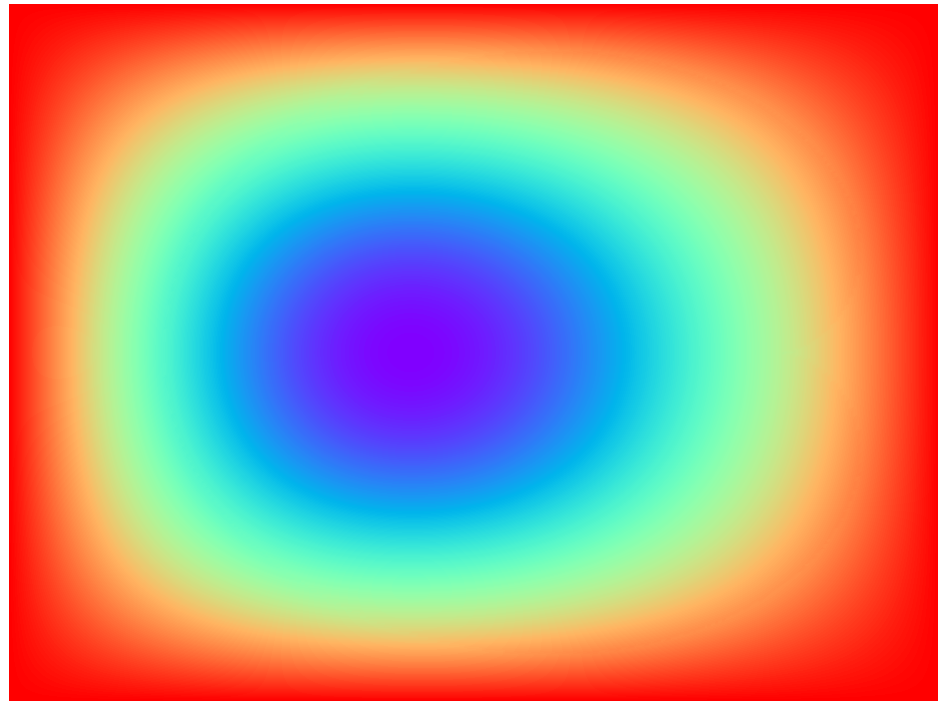} \\
		\rotatebox{90}{$r = 66$} &
		\includegraphics[width=1.2in,height=1.2in]{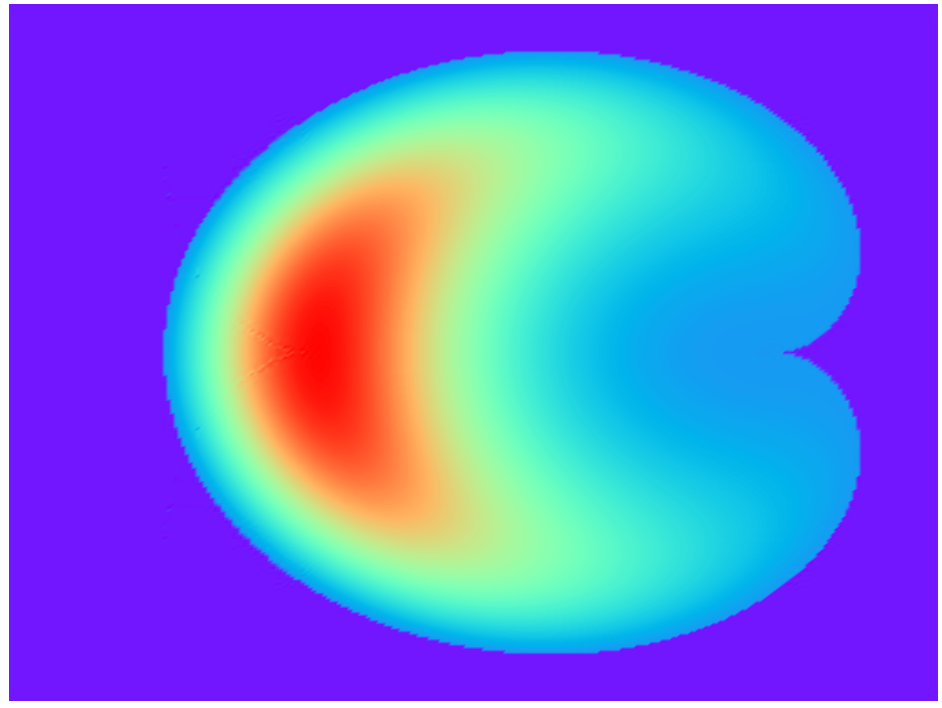} &
		\includegraphics[width=1.2in,height=1.2in]{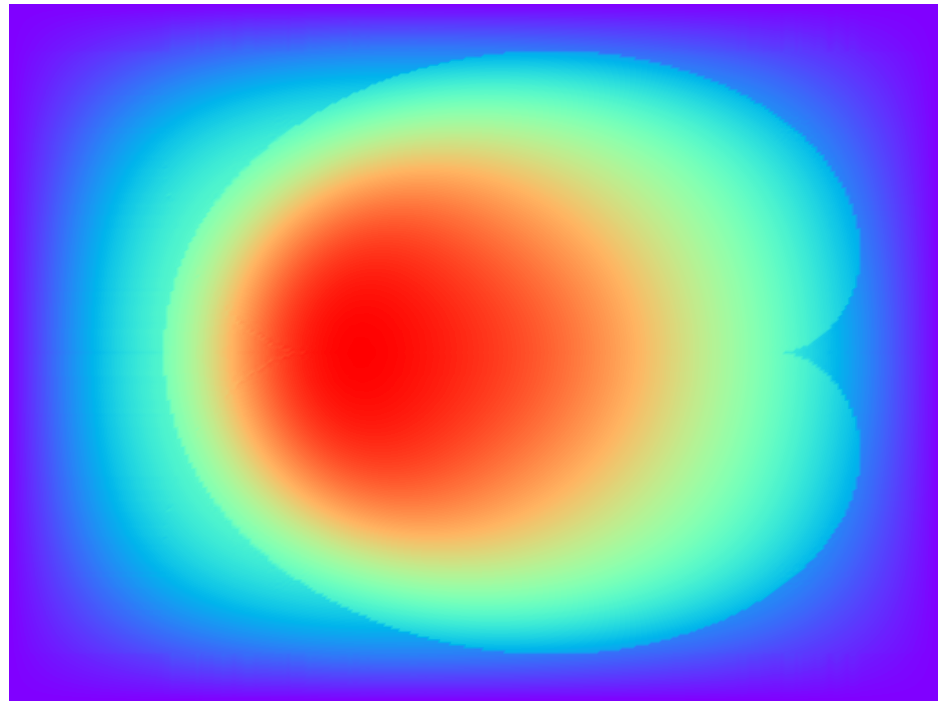} &
		\includegraphics[width=1.2in,height=1.2in]{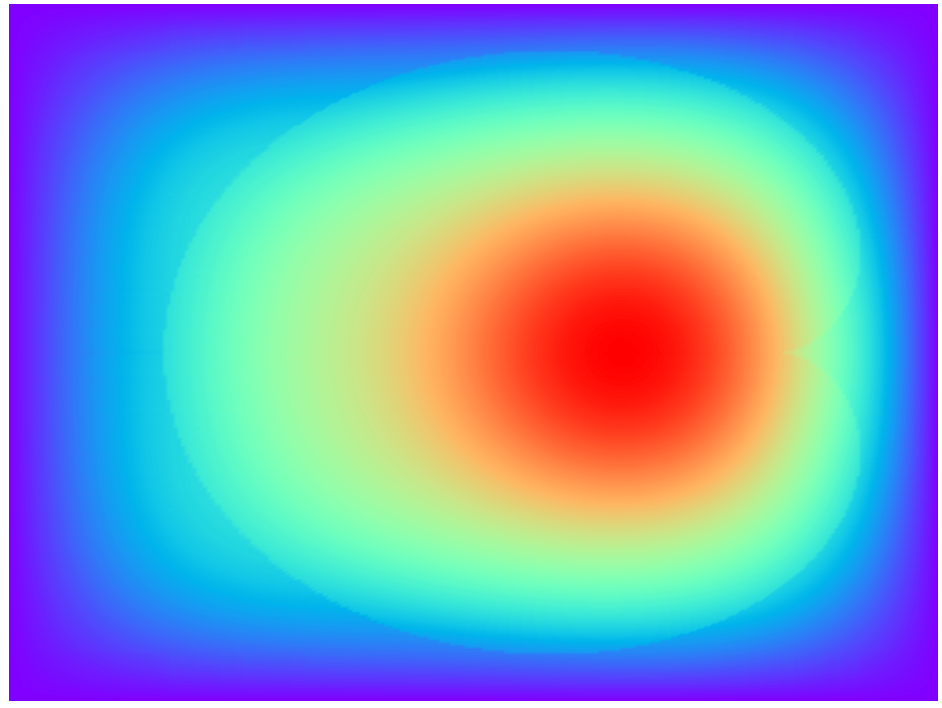} &
		\includegraphics[width=1.2in,height=1.2in]{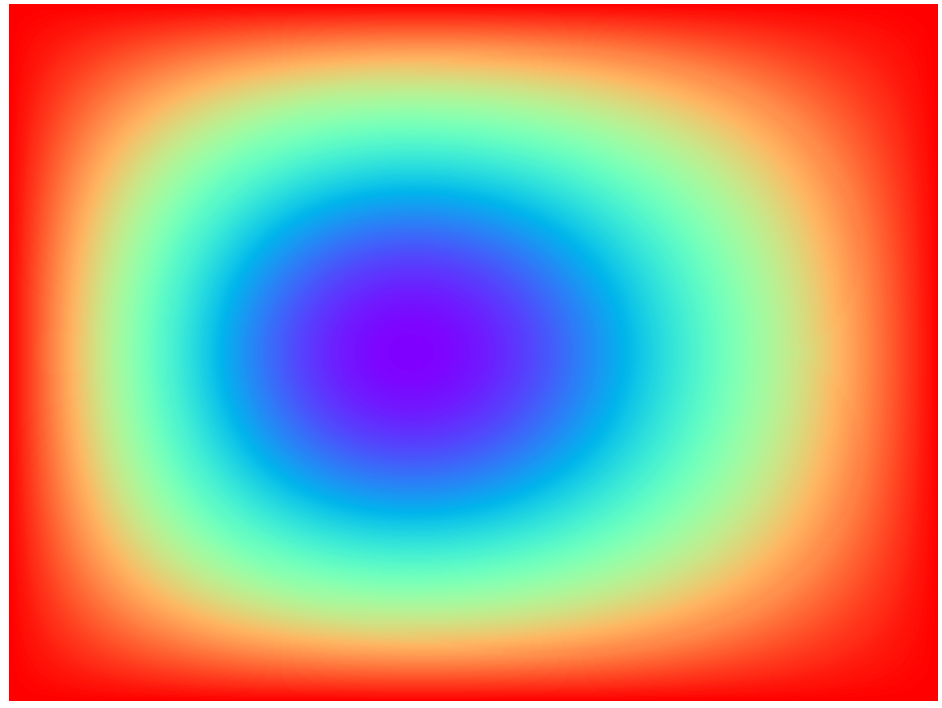} \\
		\rotatebox{90}{EI-SW22} &
		\includegraphics[width=1.2in,height=1.2in]{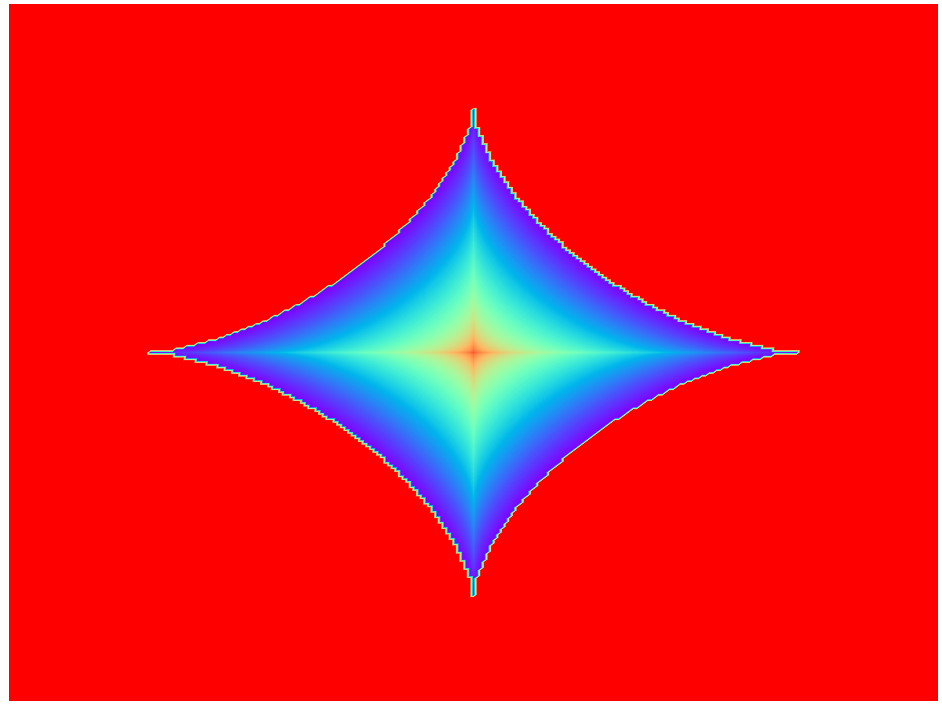} &
		\includegraphics[width=1.2in,height=1.2in]{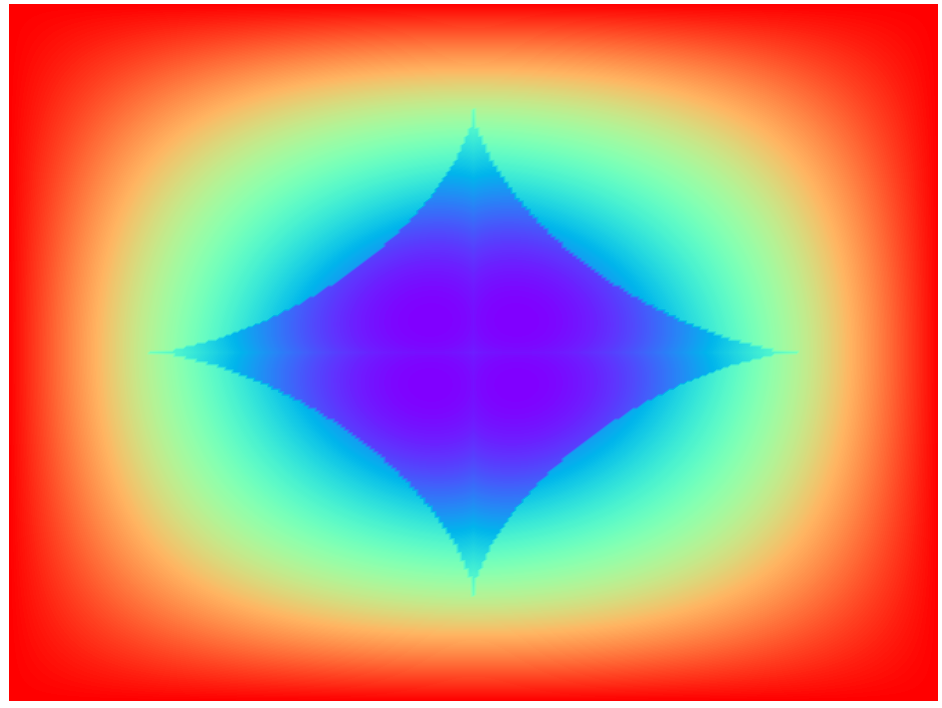} &
		\includegraphics[width=1.2in,height=1.2in]{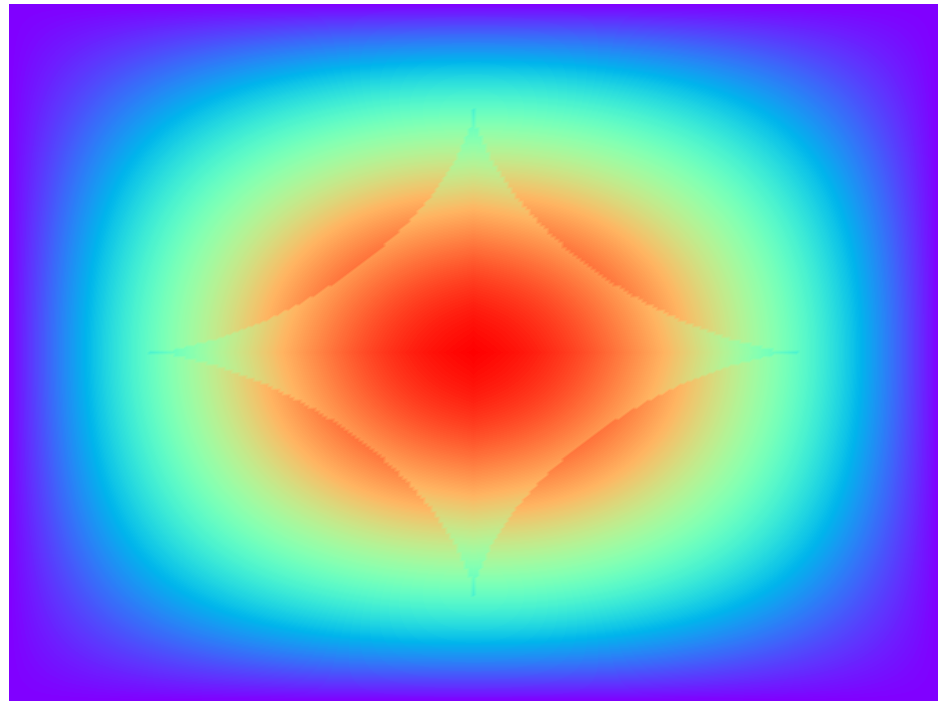} &
		\includegraphics[width=1.2in,height=1.2in]{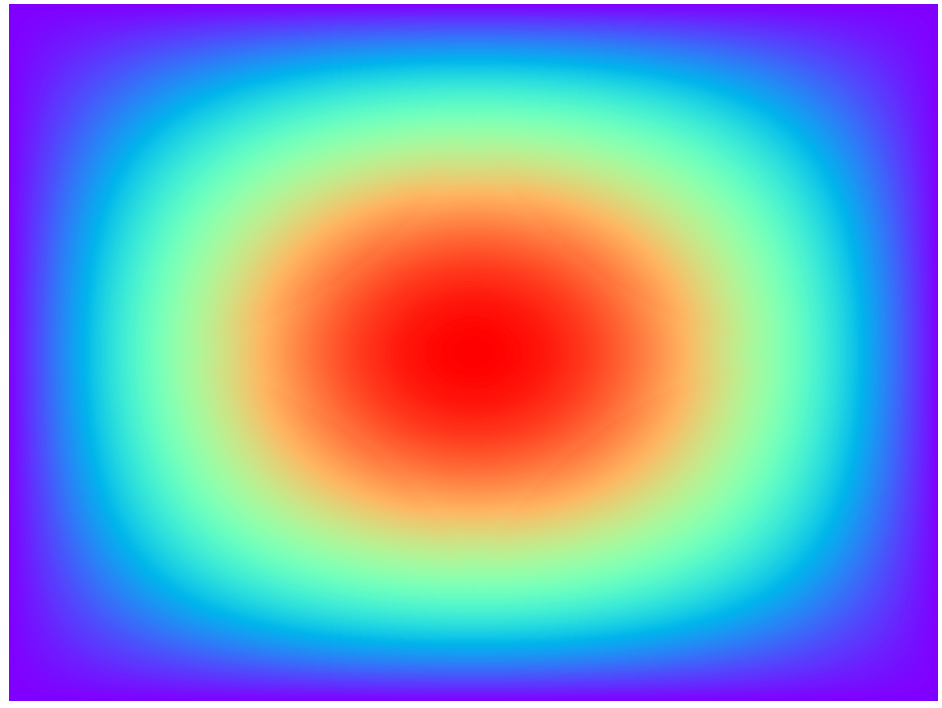} \\
		\rotatebox{90}{$r = 73$} &
		\includegraphics[width=1.2in,height=1.2in]{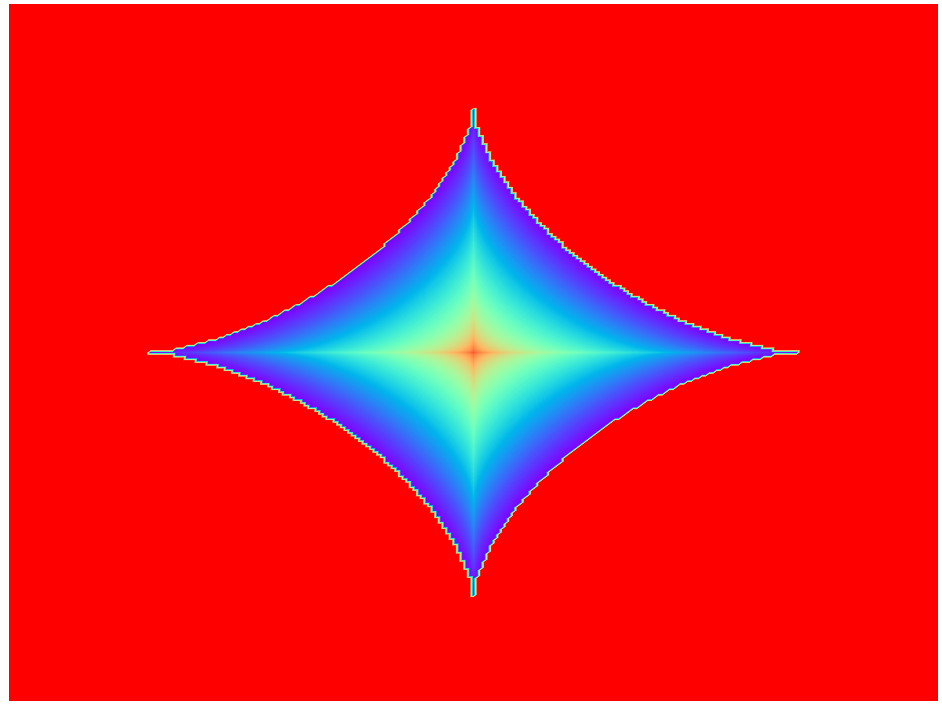} &
		\includegraphics[width=1.2in,height=1.2in]{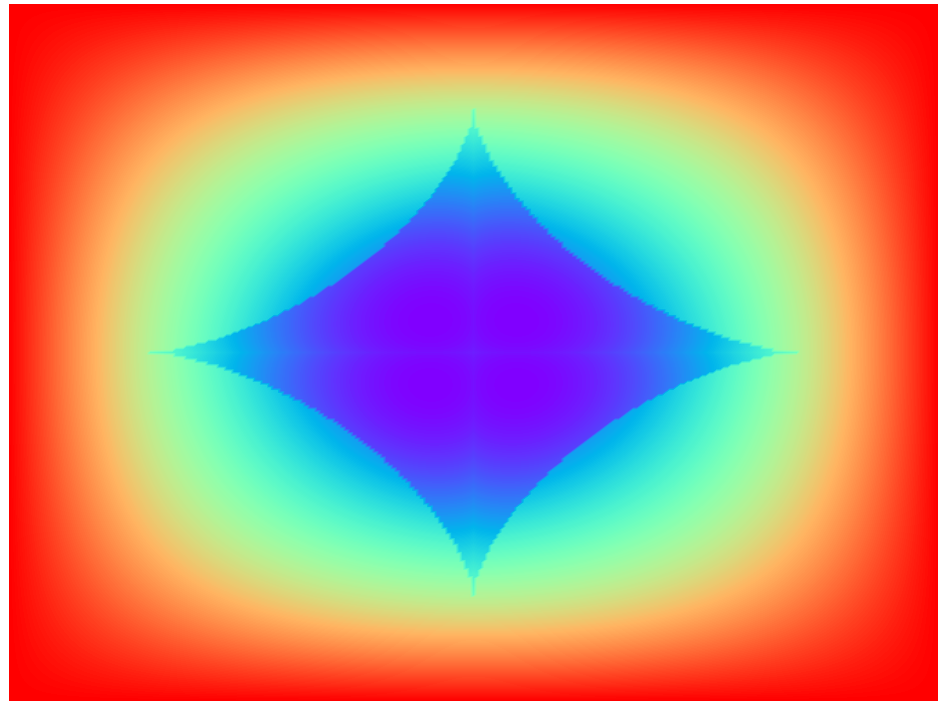} &
		\includegraphics[width=1.2in,height=1.2in]{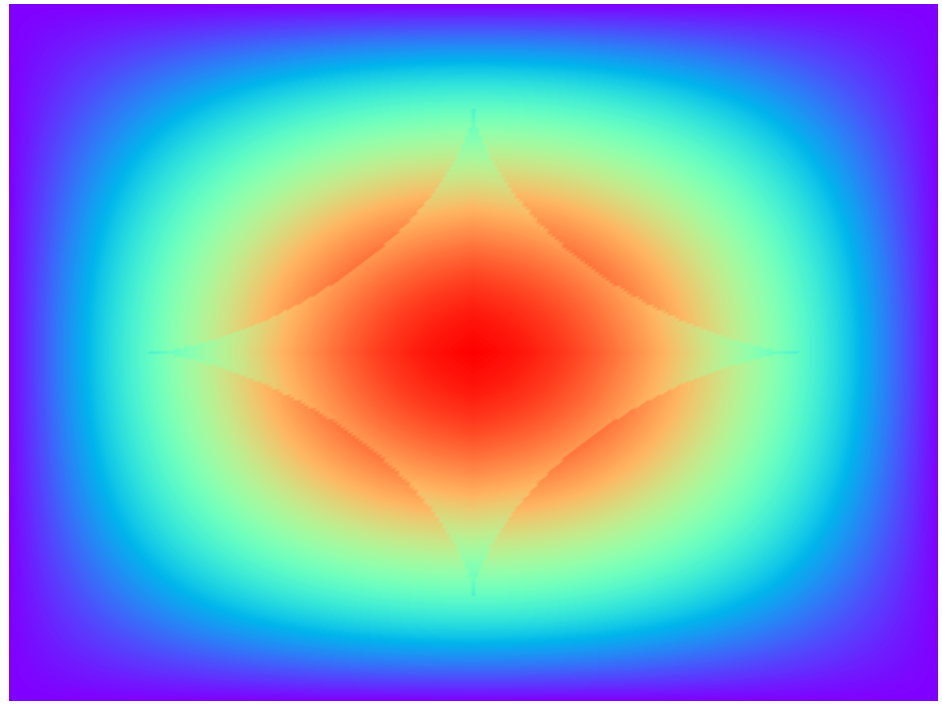} &
		\includegraphics[width=1.2in,height=1.2in]{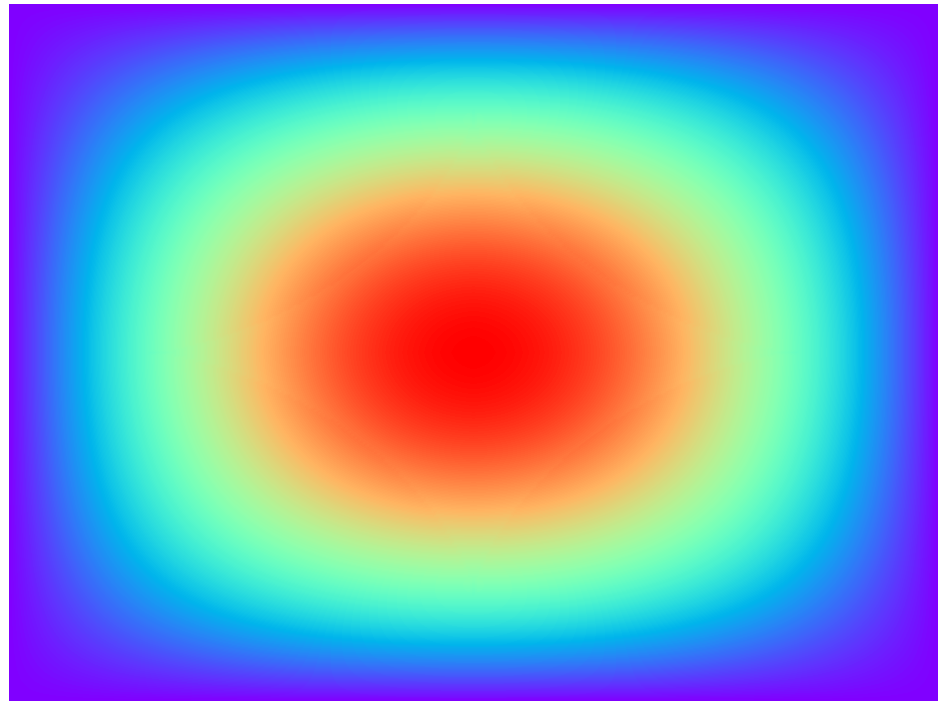} \\
	\end{tabular}
	\caption{Comparisons of the EI-SW22 \cite{phan2022exponential} 
		and our low-rank (i.e., the scheme \eqref{eq3.6}) solutions 
		for Example 3 with $(\tau,N) = (0.01,256)$.
		First two rows: 5-petals flower. Middle two rows: Cardioid. Bottom two rows: Astroid.}
	\label{fig3}
\end{figure}

\Cref{fig3} compares the numerical solutions of the EI-SW22 method and our low-rank approximation
at different moment.
Clearly, numerical solutions computed by the two methods are almost the same. 
This figure indicates that the low-rank approximation \cref{eq3.6} is reliable and accurate.
As we know, Eq.~\cref{eq1.1} is a mathematical model that 
describes the behavior of waves that lose energy over time, 
thus \cref{fig3} indeed shows that the shapes 
(i.e., 5-petals flower, cardioid and astroid) of the solutions produced 
by two numerical methods disappear over time.
Finally, these shapes change to the shape of a circle.

\section{Conclusion}
\label{sec5}
In this paper, we study low-rank solutions of strongly damped wave equations with damping term,
structural damping term and mass term.
Firstly, we employ a finite difference method to approximate the space derivatives.
A second-order matrix ordinary differential equation \cref{eq2.1} is received.
Then, we transform it to the first-order system \cref{eq2.2}.
Combining a standard Strang splitting and the DLR approximation,
our low-rank scheme \cref{eq3.6} is designed based on this first-order system \cref{eq2.2}.
Numerical results demonstrate that the proposed low-rank scheme is efficient,
reliable and accurate. 
It also implies that this scheme has a second-order convergence rate in time.
In future work, the following topics are interesting: \\
(i) Extend the DLR approximation \cite{koch2007dynamical} to stochastic differential equations or
non-local differential equations (particularly, non-local in time); \\
(ii) Develop some method (e.g., machine learning and neural networks) to search for 
a prior estimate of the best approximate rank $r$, according to the prior information of models. \\
(iii) For a blow-up problem, its solution is low-rank and sparse near the blow-up time.
Thus, it is interesting to design some method considering both the low-rank and sparsity of the solution.
We suggest a possible way:
for the blow-up problem $\dot{A}(t) = F(t, A(t)), A(t) \in \mathbb{R}^{m \times n}$, near the blow-up time,
we change to solve the following optimization problem:
\begin{equation*}
	\min\; \left\| \dot{Y}(t) - F(Y(t)) \right\| + \mathrm{sparse~term}, \qquad 
	\mathrm{s.t.}~Y(t) \in \mathcal{T}_{Y(t)} \mathcal{M}_r,
\end{equation*}
where $Y(t)$ is the rank-$r$ sparse approximation to $A(t)$, 
and $\mathcal{T}_{Y(t)} \mathcal{M}_r$ represents the tangent space of $\mathcal{M}_r$ 
at the approximation $Y(t)$.


\bibliographystyle{siamplain}
\bibliography{references}
\end{document}